\documentclass{siamltex}
\usepackage{graphicx,color}
\usepackage{amssymb,amsmath,amsbsy,mathrsfs,eucal,amsfonts}
\usepackage{algorithm}
\usepackage{algorithmic}
\def\calN{{\mathcal N}}
\def\calD{{\mathcal D}}
\def\calP{{\mathcal P}}
\def\bv{{\mathbf v}}
\def\boldf{{\mathbf f}}

\definecolor{blck}{rgb}{0,0,0}
\definecolor{darkred}{rgb}{.6,.1,0}

\newcommand\rev[1]{\textcolor{blck}{#1}}

\begin{document}

\title{Reduced Collocation Methods: Reduced Basis Methods in the Collocation Framework}
\author{Yanlai Chen \thanks{Department of
Mathematics, University of Massachusetts Dartmouth, 285 Old Westport Road, North Dartmouth, MA 02747, USA.
The research of this author was partially supported by National Science Foundation grant DMS-1216928, and by UMass Dartmouth Chancellor's Research Fund and Joseph
P. Healey Endowment Grants. (\tt{yanlai.chen@umassd.edu})} \and Sigal Gottlieb
\thanks{Department of Mathematics, University of Massachusetts Dartmouth, 285 Old Westport Road, North
Dartmouth, MA 02747, USA. The research of this author was partially supported by AFOSR grant
FA9550-09-1-0208. (\tt{sgottlieb@umassd.edu})}}
%\date%{January, 2011}
\maketitle

\begin{abstract}
In this paper, we present {\em the first} reduced basis method well-suited for the collocation framework.
Two fundamentally different algorithms are presented: the so-called Least Squares Reduced Collocation Method
(LSRCM) and Empirical Reduced Collocation Method (ERCM). This work provides a reduced basis strategy to
practitioners who \rev{prefer} a collocation, rather than Galerkin, approach. Furthermore, the empirical
reduced collocation method {\em eliminates} a potentially costly online procedure that is needed for
non-affine problems with Galerkin approach.
 Numerical results demonstrate the high efficiency and accuracy of the reduced collocation methods, which match or exceed
  that of the traditional reduced basis method in the Galerkin framework.
\end{abstract}



\section{Introduction}
\label{sec:intro}

Reduced basis methods (RBM)
\cite{Almroth_Stern_Brogan,Nagy,Noor_Peters,Prudhomme_Rovas_Veroy_Maday_Patera_Turinici, Rozza_Huynh_Patera,
CHMR_Sisc} were developed for scenarios that require a large number of numerical solutions to a parametrized
partial differential equation in a fast/real-time fashion. Examples of such situations  include simulation-based
design, parameter optimization, optimal control, multi-model/scale simulation etc. In these situations, we are
willing to expend significant computational time to pre-compute data that can be later used to compute accurate
solution in real-time.

The RBM splits the solution procedure into two parts: an offline part where the parameter dependence is
examined and a greedy algorithm is utilized to  judiciously select $N$ parameter values for pre-computation;
and an online part when the solution for any new parameter is efficiently computed based on these $N$ {\em
basis} functions.

The motivation behind the RBM is  the recognition that parameter-induced solution manifolds can be well
approximated by finite-dimensional spaces. For linear affine problems, RBM can improve efficiency by several
orders of \rev{magnitude}. For nonlinear or non-affine problems, there are remedies which allow the RBM
methods to be used efficiently
\cite{Barrault_Nguyen_Maday_Patera,Grepl_Maday_Nguyen_Patera,NguyenPateraPeraire2008}. The offline selection
of the $N$ parameter values for the pre-computed bases is enabled by a rigorous a posteriori error estimate
which guarantees the accuracy of the solution. Exponential convergence with respect to $N$ has been commonly
observed, see \cite{Rozza_Huynh_Patera, CHMR_Sisc} and the reference therein. Theoretically, a priori
convergence is confirmed for a one dimensional parametric problem \cite{Maday_Patera_Turinici_2}. More
recently, exponential convergence of the greedy algorithm for continuous and coercive problems with
parameters in any dimension has been established in \cite{BuffaMadayPateraPrudhommeTurinici2011}, and
improved in \cite{BinevCohenDahmenDevorePetrovaWojtaszczyk}.

The development and analysis of RBM has been \rev{previously} carried out in the Galerkin framework. That
is, the truth approximations (the numerical approximation from a presumably very accurate numerical scheme)
are obtained from a (Galerkin) finite element method, and the reduced basis solution is sought as a Galerkin
projection onto a low dimensional space.  However, to date RBM have not been developed, applied, or analyzed
in the context of collocation methods. While Galerkin methods are derived by requiring that the projection
of the residual onto a prescribed space is zero, collocation methods require the residual to be zero at some
pre-determined collocation points. \rev{Compared to collocation methods,
Galerkin methods have a weaker regularity requirement on the solution. For example, for second-order problems,
collocation methods require the solution to be at least $H^2$ over the domain $\Omega$, while
Galerkin methods only require  solutions in $H^1(\Omega)$, due to the adoption of the weak formulation.
Unlike collocation methods, Galerkin methods do not require a tensorial grid and handle curved boundaries with ease.
On the other hand,} collocation methods are particularly attractive for their ease of
implementation, particularly for
 time-dependent nonlinear problems \cite{TrefethenSpecBook,ShenTangBook,HesthavenGottlieb2007}.
%\rev{Its limits in applicability compared to a Galerkin approach mainly include  stronger regularity requirement in the field variables and the typical domains (or subdomains) being tensorial. On the other hand, Galerkin approaches can, for an example, consider solutions in $H^1(\Omega)$ for
%second-order problems due to the adoption of the weak formulation. Moreover, it deal with curved boundary with ease.}

In this paper, we develop the RBM idea for collocation methods. Given a highly accurate
collocation method that is used as the truth solver for the parametric problem, we wish to
study the performance of the system under variation of certain parameters
using a collocation-based RBM. That is, the new method uses
collocation for both the truth solver and the online reduced solver.

The paper is organized as follows. In Section \ref{sec:alg}, we present two approaches to collocation RBM.
The first one utilizes a least squares approach. The second one relies on a projection of the fine
collocation grid problem onto a (carefully-chosen) coarse collocation grid.
Theoretical analysis and discussions on the offline-online decomposition are provided in Section
\ref{sec:analysis}. Numerical results are shown in Section \ref{sec:numerical}. Finally, some
concluding remarks and future directions are  in Section \ref{sec:conclude}.

\section{The Algorithms}
\label{sec:alg}

We begin with  a parametrized partial differential equation of the form
\begin{equation}
\mathbb{L} (\mu)\, u_\mu (x) = f(x; \mu), \qquad x \in \Omega \,\,\,\rev{\subset \mathbb{R}^n} \label{eq:pde}
\end{equation}
with appropriate boundary conditions. We are interested in the solutions of the differential equation over a
range of parameter values $\mu$, where $\mu = (\mu_1, \dots, \mu_d) \in \calD$, a prescribed
\rev{$d$-}dimensional real parameter domain. The parameters can be, for example, heat conductivity, wave
speed, angular frequency, or geometrical configurations etc.

In this work, we assume that the operator is {\em linear} and  {\em affine} with respect to functions of $\mu$.
That is, $\mathbb{L} (\mu)$ can be written
as a linear combination of parameter-dependent coefficients and parameter-independent operators:
\begin{equation}
\mathbb{L} (\mu) = \sum_{q = 1}^{Q_a} a^{\mathbb{L}}_q(\mu) \mathbb{L}_q. \label{eq:affineL}
\end{equation}
We make a similar assumption for $f$:
\begin{equation}
f(x; \mu) = \sum_{q = 1}^{Q_f} a^f_q(\mu) f_q (x). \label{eq:affinef}
\end{equation}
In the Galerkin framework, these are common assumptions in the reduced basis literature
\cite{Rozza_Huynh_Patera}. There are remedies available when the parameter-dependence is not affine
\cite{Barrault_Nguyen_Maday_Patera,Grepl_Maday_Nguyen_Patera,NguyenPateraPeraire2008}.

For any value of the parameter $\mu$, we can approximate the solution to this equation
using a collocation approach: we define a discrete differentiation operator $\mathbb{L}_\calN (\mu)$
so that the approximate solution $u^\calN_\mu$ satisfies the equation
\begin{equation}
\mathbb{L}_\calN (\mu)\, u^\calN_\mu (x_j) = f(x_j; \mu), \label{eq:discpde}
\end{equation}
{\em exactly} on  a given set of $\calN$ collocation points $C^\calN = \{x_j\}_{j=1}^\calN$, \rev{usually
taken as a tensor product of $\calN_x$ collocation points for each dimension. Obviously, for $\Omega \subset \mathbb{R}^n$ we have $\calN =
\calN_x^n$.} We assume that the scheme \eqref{eq:discpde} produces highly accurate numerical solutions
$u^\calN_\mu$ to the problem \eqref{eq:pde}. We refer to the solution $u^\calN_\mu$ as the ``truth
approximation''.

Although solving   \eqref{eq:discpde} gives highly accurate approximations, it is prohibitively expensive and
time-consuming to repeat for a very large number of parameter values  $\mu$. The reduced basis method allows
for highly accurate solutions to be computed quickly and efficiently when needed (the ``online''
computation) based on a set of possibly expensive offline computations. The idea of the reduced basis method
is that we first {\em pre-compute} the truth approximations for a set of \rev{$N << \calN$} well-chosen parameter values
$\{ \mu^1,\, \mu^2, \, \dots,\, \mu^N \}$ by solving \eqref{eq:discpde} with the corresponding parameter value.
Then when the solution for any parameter value $\mu^*$ in the (prescribed) parameter domain $\calD$ is
needed, instead of solving for the (usually expensive) truth approximation $u^\calN_{\mu^*}$, we combine
$u^\calN_{\mu^1},\, u^\calN_{\mu^2}, \, \dots,\, u^\calN_{\mu^N}$ in some way to produce a surrogate
solution $u^{(N)}_{\mu^*} $:
\[ u^{(N)}_{\mu^*} = \sum_{j=1}^{N} c_j(\mu^*) u^\calN_{\mu^j}. \]

Thus, the design of the reduced basis method requires two components:
\begin{enumerate}
\item Offline: how to select the pre-computed basis.
\item Online: how to combine the pre-computed basis functions to produce the surrogate solution.
\end{enumerate}
In the following sections, we describe two variants of the reduced collocation algorithm. We first explain
our approaches for the online computation of the surrogate solution from the pre-computed reduced basis in
Section \ref{sec:chooseprojection}, and then the related question of the selection of the reduced basis in
Section \ref{sec:precomp}.

\subsection{Online algorithms}
\label{sec:chooseprojection}

For the surrogate solution $u^{(N)}_{\mu^*} $ to approximate the truth approximation $u^\calN_{\mu^*}$ reasonably
well, we require that $u^{(N)}_{\mu^*} $ provides, in some sense, a good approximation to the solution of the
discretized differential equation
\[ \mathbb{L}_\calN(\mu^*)  \left(\sum_{j=1}^{N} c_j u^\calN_{\mu^j}\right)  \approx f(x; \mu^*).\]
By exploiting the linearity of the operator we observe that our task is to find coefficients $c_j(\mu^*)$ so
that the residual
\[ \sum_{j=1}^{N} c_j(\mu^*) \mathbb{L}_\calN(\mu^*) u^\calN_{\mu^j}  - f(x; \mu^*) \]
is small.
%(With a slight abuse  of notation, we use $u_\mu^\calN$ to denote both the function and the vector of
%the function values at the grid points.)

In the Galerkin framework \cite{Rozza_Huynh_Patera,CHMR_Sisc}, the coefficients \rev{are} found by requiring
that the $L^2$-projection of this residual onto the reduced space is zero. For the collocation case, the
system of equations we \rev{wish} to solve is
\begin{equation}
\sum_{j=1}^{N} c_j(\mu^*) \mathbb{L}_\calN(\mu^*) u^\calN_{\mu^j}(x_k)  =  f(x_k; \mu^*) \; \; \; \; k=1,  .
. ., \calN.
\end{equation}
However, this system is over-determined: we have  only $N$ unknowns, but $\calN >> N$ equations. To
approximate the solution to this system, our task is to identify an appropriate \rev{operator} $\calP$ such that
the following holds
\begin{equation}
\calP\left(\sum_{j=1}^{N} c_j(\mu^*) \mathbb{L}_\calN(\mu^*) u^\calN_{\mu^j}\right)  = \calP\left(f(x;
\mu^*)\right). \label{eq:rbmpde}
\end{equation}
By considering two different ways to choose the operator $\calP$ in Equation \eqref{eq:rbmpde},
we propose two approaches for finding the coefficients of the reduced basis solutions. These
two approaches are the least squares approach and the reduced collocation method.

\noindent{\bf Least-squares approach.} Our first approach is a very standard approach to approximating the
solution to an over-determined system. We  determine the coefficients by satisfying the equation
\eqref{eq:rbmpde} in a least squares sense. Given $\{u_{\mu^1}^\calN$, $u_{\mu^2}^\calN$, $\dots$,
$u_{\mu^N}^\calN\}$, we define, for any $\mu^*$, an $\calN \times N$ matrix
\[
\mathbb{A}_N (\mu^*) = \left(\mathbb{L}_\calN(\mu^*)\, u_{\mu^1}^\calN, \mathbb{L}_\calN(\mu^*)\,
u_{\mu^2}^\calN, \, \dots, \, \mathbb{L}_\calN(\mu^*)\, u_{\mu^N}^\calN\right),
\]
and vector of length $\calN$
\[ \boldf^\calN_j = f(x_j ; \mu^*) \; \; \; \; x_j \in C^\calN, \]
and solve \rev{the least squares problem}
\begin{equation}
\mathbb{A}_N^T(\mu^*)\,\mathbb{A}_N(\mu^*) \,\vec{c} = \mathbb{A}_N^T(\mu^*)\, \boldf^\calN  \label{eq:reducedprobLS}
\end{equation}
to obtain $\vec{c} = (c_1(\mu^*), c_2(\mu^*), \dots, c_N(\mu^*))^T$.

\noindent{\bf Reduced Collocation approach.}
A  more natural approach from the collocation point-of-view  is to determine the coefficients $c$ by
enforcing \eqref{eq:rbmpde} at a reduced set of collocation points $C_R^N$. In other words, we solve
\begin{equation}
\sum_{j=1}^{N} c_j(\mu^*) \mathbb{\rev{I}}_\calN^N \left(\mathbb{L}_\calN(\mu^*) u^\calN_{\mu^j}\right)  = f(x;
\mu^*). \quad {\rm for} \quad x \in C_R^N\label{eq:rbmpderc}
\end{equation}
where $\mathbb{\rev{I}}_\calN^N$ is the \rev{interpolation} operator \rev{for  functions defined in} the $\calN$-dimensional space corresponding to
the fine-domain collocation points $C^\calN$ \rev{at} the smaller %$N$-dimensional space associated with the
reduced set of collocation points
$C^N_R$.
In other words, we  define the $N$ vectors of length $N$ by their elements
\[ \left(\bv^{j}_{\mu^*}\right)_k = \left. \left( \mathbb{L}_\calN(\mu^*) u^\calN_{\mu^j}\right)\right|_{x=x_k} \; \; \;
\mbox{and} \; \; \; \; \left(\boldf^N\right)_k = f(x_k) \; \; \; \mbox{for $ x_k \in  C_R^N$}\] and solve
the $N \times N $ system of equations
\begin{equation}
\sum_{j=1}^{N} c_j(\mu^*) \bv^{j}_{\mu^*} = \boldf^N.
\end{equation}
The choice of reduced collocation points $C_R^N$
%is not limited to the Chebyshev grid, and in principle
can be any set of $N$ points in the computational domain. Later we will demonstrate how this set of points can
be determined, together with the choice of basis functions, through the greedy algorithm  (Algorithm
\ref{alg:RCgreedy}).
Although the coefficients are computed based on collocation on  a coarser mesh,
the quality of the reduced solution is not degraded since the differentiations are performed first,
by the highly accurate operator $\mathbb{L}_\calN(\mu^*)$  whose accuracy is dependent on $\calN$.
This differentiation is then followed by \rev{an interpolation at} the set of $N$ points.

Once the coefficients $\{c_j(\mu^*)\}$ are determined, whether by the least squares approach or the reduced
collocation approach, we define the reduced basis solution
\[ u^{(N)}_{\mu^*} = \sum_{j=1}^N\, c_j(\mu^*) u_{\mu^j}^\calN. \]
In both \rev{the least squares and reduced collocation cases}, the coefficients are determined by solving an $N \times N$ system. Furthermore, due to the
affine assumption on the operator, the online cost of assembling the system is also independent of  $\calN$
(as will be seen in Section \ref{sec:analysis}). Thus, the online component requires only modest
computational cost because $N$ is not large.

\subsection{The pre-computation stage}
\label{sec:precomp} Appropriate selection of the basis functions is a major determinant of how well the
reduced basis method will work. The pre-computation and selection of basis solutions may be expensive and
time-consuming, but this cost is acceptable because it is offline and done once-for-all. Once the reduced
basis solutions are computed and selected, the online component can proceed efficiently, as described above.

In this section we describe algorithms for choosing the reduced basis set $\{ u_{\mu^1}^\calN, \dots,
u_{\mu^N}^\calN \}$. The selection of the
reduced basis is performed in order to enable us to certify the accuracy of  the reduced solution. The
critical piece of information is that  given a pre-computed reduced basis set  $\{ u_{\mu^1}^\calN, \dots,
u_{\mu^i}^\calN \}$ we can compute  an upper bound $\Delta_i(\mu)$  for the error of the reduced solution
$u_\mu^{(i)}$ for any parameter $\mu$. This upper bound is given by
\begin{equation}
\Delta_i(\mu) = \frac{\lVert \rev{\boldf^\calN - \mathbb{L}_\calN (\mu)
u_\mu^{(i)}}\rVert_{\ell^2}}{\rev{\sqrt{\beta_{LB}(\mu)}}}.
\end{equation}
where $\beta_{LB}(\mu)$ is the lower bound for  the smallest eigenvalue of $\mathbb{L}_\calN (\mu)^T
\mathbb{L}_\calN (\mu)$. This upper bound is enabled by the {\em a posteriori} error estimate which will be
proved in Section \ref{sec:apost_greedy}.
In the following, we present  the greedy algorithms used for the selection of the pre-computed basis
for the least squares and the reduced collocation approaches.

\subsubsection{Least Squares Reduced Collocation Method (LSRCM)}
\label{sec:LSRCM}

The idea behind the greedy algorithm is to discretize the parameter space, and scan the discrete parameter
space to select the best reduced solution space. To do this, we first randomly select one parameter and call
it $\mu^1$, and compute the associated highly accurate solution $u^\calN_{\mu^1}$. Next, we scan the entire
discrete parameter space and for each parameter in this space compute its least squares reduced basis
approximation $u^{(1)}_{\mu}$. We now compute the error \rev{estimator} $\Delta_1(\mu)$. The next parameter
value we select, $\mu^2$, is the one corresponding to the largest error estimator.
 \rev{We then compute the  highly accurate solution $u_{\mu^2}^\calN$, and thus
have  a new basis set consisting of two elements} $\{ u_{\mu^1}^\calN, u_{\mu^2}^\calN\}$

\begin{algorithm}[htp]
  \caption{Least Squares Reduced Collocation Method (LSRCM)\rev{: Offline Procedure}}\label{alg:LSgreedy}
  \begin{algorithmic}%[1]
\STATE {\bf 1.} \rev{Discretize the parameter domain $\calD$ by $\Xi$, and denote the center of $\calD$ by $\mu^c$.}
\medskip
\STATE {\bf 2.} Randomly select $\mu^1$ and solve $\mathbb{L}_\calN (\mu^1)\, u^\calN_{\mu^1}
(x) = f(x; \mu^1)$ for $x \in C^\calN$.
\medskip
\STATE {\bf 3.} For $i = 2, \dots, N$ do%\IF{sth} \STATE do this; \ENDIF
\begin{itemize}
\item [{\bf 1).}] Form $\mathbb{A}_{i-1} = \left(\mathbb{L}_\calN\, u_{\mu^1}^\calN, \mathbb{L}_\calN\, u_{\mu^2}^\calN, \, \dots, \,
\mathbb{L}_\calN\, u_{\mu^{i-1}}^\calN\right)$.
\item [{\bf 2).}] For all $\mu \in \Xi$, solve
$\mathbb{A}_{i-1}^T\,\mathbb{A}_{i-1} \,\vec{c} = \mathbb{A}_{i-1}^T\,\boldf^\calN$ to obtain $u^{(i-1)}_\mu
= \sum_{j = 1}^{i - 1} c_j u^\calN_{\mu^j}$.
\item [{\bf 3).}] For all $\mu \in \Xi$, calculate $\Delta_{i-1}(\mu)$.
\item [{\bf 4).}] Set $\mu^i = argmax_{\mu}\,\,\Delta_{i-1}(\mu)$.
\item [{\bf 5).}] Solve $\mathbb{L}_\calN (\mu^i)\, u^\calN_{\mu^i}
(x) = f(x; \mu^i)$ for $x \in C^\calN$.
\end{itemize}
\medskip
\STATE {\bf 4.} \rev{Apply a modified Gram-Schmidt transformation, with inner product defined by $(u,v) \equiv \left(\mathbb{L}_\calN (\mu^c)u, \mathbb{L}_\calN (\mu^c) v\right)_{L^2(\Omega)}$, on the basis $\left\{u^\calN_{\mu^1}, u^\calN_{\mu^2}, \dots, u^\calN_{\mu^N}\right\}$  to obtain a more stable basis $\left\{\xi_1^\calN, \xi_2^\calN, \dots, \xi_N^\calN\right\}$ for the least squares reduced collocation method.}
\medskip
  \end{algorithmic}
\end{algorithm}

This process is repeated until the \rev{maximum of the error estimators} is sufficiently small.
 At every step we select the parameter which is approximated most badly by the current solution space,
\rev{with the goal being} that in this way we select a solution space that will approximate any parameter
reasonably well. The detailed algorithm is provided in Algorithm \ref{alg:LSgreedy}. \rev{To ensure the reduced system is well-conditioned, we apply the modified Gram-Schmidt transformation with weighted inner product.}

%\clearpage

\subsubsection{Empirical Reduced Collocation Method (ERCM)}
\label{sec:ERCM}

The least squares approach above can not be immediately adapted to the collocation case because collocation
requires the same number of collocation points as basis functions. Thus we face \rev{an additional} problem
of having to choose an appropriate set of collocation points $C_R^N$ at which to enforce the PDE.
In fact, the choice of the reduced set of collocation points is crucial for the  accuracy of the algorithm.
For example, as we will show in the numerical example in Section \ref{sec:numerical},
naively using the coarse Chebyshev grid does not yield an accurate result.  In the following, we
propose the {\em Empirical Reduced Collocation Method} for choosing the basis functions and reduced
collocation points.

\begin{algorithm}[h!]
  \caption{Empirical Reduced Collocation Method (ERCM)\rev{: Offline Procedure}}\label{alg:RCgreedy}
  \begin{algorithmic}%[1]
\STATE {\bf 1.} Randomly select $\mu^1$ and solve $\mathbb{L}_\calN (\mu^1)\, u^\calN_{\mu^1}
(x) = f(x; \mu^1)$,
and let $x^1 = argmax_{x \in X} \,\, \left|u^\calN_{\mu^1}(x)\right|,\quad \xi^\calN_1 =
\frac{u^\calN_{\mu^1}}{u^\calN_{\mu^1}(x^1)}$.
\medskip
\STATE {\bf 2.} For $i = 2, \dots, N$ do%\IF{sth} \STATE do this; \ENDIF
\begin{itemize}
\item [{\bf 1).}] Let $C_R^{i-1} = \left\{x^1,\dots,x^{i-1}\right\}$.
\item [{\bf 2).}] For all $\mu \in \Xi$, solve \ $\sum_{j=1}^{i-1} c_j \mathbb{P}_\calN^N
\left(\mathbb{L}_\calN(\mu) u^\calN_{\mu^j}\right)  = f(x; \mu) \,\, {\rm for} \,\, x \in C_R^{i-1}$ to
obtain $u^{(i-1)}_\mu = \sum_{j = 1}^{i - 1} c_j u^\calN_{\mu^j}$.
\item [{\bf 3).}] For all $\mu \in \Xi$, calculate $\Delta_{i-1}(\mu)$.
\item [{\bf 4).}] Set $\mu^i = argmax_{\mu \in \Xi}\,\,\Delta_{i-1}(\mu)$.
\item [{\bf 5).}] Solve $\mathbb{L}_\calN (\mu^i)\, u^\calN_{\mu^i}
(x) = f(x; \mu^i)$.
\item [{\bf 6).}] Find $\alpha_1, \dots, \alpha_{i-1}$ such that, if we define $\xi^\calN_i = u^\calN_{\mu^i} - \sum_{j=1}^{i-1}
\alpha_j\,\xi^\calN_j$, we have $\xi^\calN_i(x^j) = 0$ for $j = 1, \dots, i-1$.
%\item [{\bf 2).}] Set $\xi^\calN_i(x) = 0$ for $x$ in the neighborhoods of $x^1, \dots, x^{i-1}$.
\item [{\bf 7).}] Set $x^i = argmax_x \,\, \left|\xi^\calN_{i}\right|$ and $\xi^\calN_i = \frac{\xi^\calN_i}{\xi^\calN_i(x^i)}$.
\item[\rev{\bf 8).}] \rev{Apply modified Gram-Schmidt transformation on $\left\{\xi_1^\calN, \dots, \xi_i^\calN\right\}$.}
\end{itemize}
\medskip
\STATE {\bf 3.} Set the reduced set of collocation points $C_R^N = \left\{x^1, x^2, \dots, x^N\right\}$ and
%use the $\xi^\calN_1, \dots, \xi^\calN_N$ as the set of bases for the reduced method.
use the set
\rev{$\left\{\xi_1^\calN, \xi_2^\calN, \dots, \xi_i^\calN\right\}$ as the basis for the reduced collocation method.}
%$\left\{u^\calN_{\mu^1}, u^\calN_{\mu^2}, \dots, u^\calN_{\mu^N}\right\}$ as the basis for the reduced collocation method.
\medskip
  \end{algorithmic}
\end{algorithm}
The idea behind the empirical reduced collocation method is similar to the greedy algorithm used quite often
by \rev{the} reduced basis method and it has the same structure as the Empirical Interpolation Method
\cite{Barrault_Nguyen_Maday_Patera,Grepl_Maday_Nguyen_Patera,NguyenPateraPeraire2008}. We build the set of collocation points
hierarchically with the each point chosen from the set of candidate points $X$ \rev{(taken to be the set of
fine collocation grid $C^\calN$)}. \rev{We begin by picking a parameter $\mu^1$ randomly and computing
the corresponding  basis function $u^\calN_{\mu^1}$, and selecting the collocation point
$x_1$ at which the absolute value of the basis function attains its maximum.
 (We note that  it is also possible to choose the
collocation point which maximizes one of the partial derivative of the basis function
however, this choice did not perform well in numerical tests.)
Now we can say we have a  set of basis functions $ \{ u^\calN_{\mu^j} \}_{j=1}^{i-1}$ and
a set of collocation points  $ \{ x_j\}_{j=1}^{i-1} $. We compute the reduced basis solution
$u^{(i-1)}_{\mu}$ for all $\mu$ in the discretized parameter domain, and the associated
error estimator $\Delta_{i-1}(\mu)$.
To get the next basis function, we find the parameter value $\mu^i$ at which the
error estimator is maximized, and we compute the highly accurate solution $u^\calN_{\mu^i}$.
To ensure well-conditioning of the process, we orthonormalize the basis functions  in the sense that if we
define $B^N_{ij} = \xi^\calN_j(x^i)$, then the matrix $B^N$ is lower triangular with unit diagonal.
We now obtain the set of orthonormalized basis functions $\{\xi^\calN_1, \dots, \xi^\calN_N\}$.
Finally, the $i$th collocation point is chosen to be that at which the absolute value of the
basis function $\xi_i^\calN$ is maximized.
Repeatedly following this procedure, given in Algorithm \ref{alg:RCgreedy},
we obtain the set of orthonormalized basis functions $\{\xi^\calN_1, \dots, \xi^\calN_N\}$ and the {\em
reduced} set of collocation points $C_R^N = \left\{x^1, x^2, \dots, x^N\right\}$ that will be used to find
the surrogate solution.}

\rev{\noindent {\bf Remark} The choice of collocation points described above is different from the ``best point''
and ``hierarchical point'' approximations described in \cite{NguyenPateraPeraire2008}. In our approach, we
used the rather ad-hoc -- and inexpensive --  approach of choosing a collocation point which maximizes the
corresponding basis function. The ``best point'' and ``hierarchical point'' approaches choose the interpolation
points by minimizing, in some sense,  the difference between the interpolation and projection coefficients.
However, we compared these approaches in numerical tests based on the problems considered in Section 4, and
the sophisticated  ``best point'' and ``hierarchical point'' approaches performed no better
 than the simple algorithm above
in terms of size of errors and rate of convergence of the reduced collocation solution to the truth approximation.}

\section{Analysis of the Reduced Collocation Method}
\label{sec:analysis}

In this section, we provide some analysis of the proposed algorithms and some details for the
offline-online decomposition that is crucial to the traditional tremendous speedup of reduced basis method.

\subsection{A Posteriori Error Estimate}
\label{sec:apost_greedy}

The essential ingredient of the accuracy of the reduced collocation method is the
upper bound which is used for error estimation. In this section,
we state and prove the theorem relating to this error estimator.

Before we state our theorem, we must assume that we
have a lower bound $\beta_{LB}(\mu)$ for the smallest eigenvalue of $\mathbb{L}_\calN (\mu)^T
\mathbb{L}_\calN (\mu)$,
\begin{equation}
\label{eq:betamu} \beta(\mu) = \min_{v}\frac{v^T\, \mathbb{L}_\calN (\mu)^T \mathbb{L}_\calN (\mu) \,v}{v^T
v}.
\end{equation}
%and the equivalence of $L^2-$norm of a function $u$ and the discrete $\ell^2-$norm of the corresponding
%vector of function values at the collocation points
%\begin{equation}
%C_1 \lVert u \rVert_{\ell^2} \le \lVert u \rVert_{0} \le C_2 \lVert u \rVert_{\ell^2}.
%\end{equation}

\begin{theorem}
For any $\mu$, suppose $u_\mu^\calN$ is the truth approximation solving \eqref{eq:discpde} and $u_\mu^{(N)}$ is
the reduced basis solution solving \eqref{eq:rbmpderc} or \eqref{eq:reducedprobLS}, we define
\begin{equation}
\Delta_N(\mu) = \frac{\lVert\rev{\boldf^\calN - \mathbb{L}_\calN (\mu)
u_\mu^{(N)}}\rVert_{\ell^2}}{\rev{\sqrt{\beta_{LB}(\mu)}}}. \label{eq:errorestimator}
\end{equation}
Then we have $\lVert u_\mu^\calN - u_\mu^{(N)}\rVert_{\ell^2} \le \Delta_N(\mu)$.
\end{theorem}
\begin{proof}
We have the following error equation on the $\calN$-dependent fine domain collocation grid thanks to the
equation satisfied by the truth approximation \eqref{eq:discpde}:
\[
\mathbb{L}_\calN (\mu) \left(u_\mu^\calN - u_\mu^N\right) = f - \mathbb{L}_\calN (\mu) u_\mu^N.
\]
%Multiplying both sides by the inverse of $\mathbb{L}_\calN (\mu)^T \mathbb{L}_\calN (\mu)$ and
Taking the discrete $\ell^2$-norm and using basic properties of eigenvalues gives
\[ \lVert u_\mu^\calN - u_\mu^N\rVert_{\ell^2} \le \frac{\lVert f - \mathbb{L}_\calN
(\mu) u_\mu^N\rVert_{\ell^2}}{\rev{\sqrt{\beta_{LB}(\mu)}}}.\]
\end{proof}

This a posteriori error estimate is used repeatedly in the
greedy algorithm to determine the reduced basis set $\left\{u^\calN_{\mu^1}, u^\calN_{\mu^2}, \dots,
u^\calN_{\mu^N}\right\}$ . In addition, the {\em a posteriori} error estimate
also serves the role of certifying the accuracy of the reduced solution: given a tolerance
$\epsilon_{\mbox{tol}}$, it is trivial to modify the algorithms so that they will find an appropriate number
$N$ and a corresponding set $\left\{u^\calN_{\mu^1}, u^\calN_{\mu^2}, \dots, u^\calN_{\mu^N}\right\}$ such
that the resulting reduced solver will have error below $\epsilon_{\mbox{tol}}$ for $\mu \in \Xi$.
While this is not enough to guarantee accuracy for any $\mu \in \calD$, it suggests that if $\Xi$ is a discretization
that represents $\calD$ well, the reduced basis method will work well for any $\mu \in \calD$.

\subsection{Offline-Online decomposition}
\label{sec:OffOn}

As is well-known \cite{Rozza_Huynh_Patera}, the tremendous speedup of the reduced basis method comes from
the decomposition of the computation into two-stages, called offline and online stages. The offline stage is
done once for all and is $\calN$-dependent (thus expensive). The online stage should be independent of
$\calN$ thus economical and can be afforded for every new value of the parameter $\mu$ in the prescribed
domain $\calD$.

Thus the key to the efficiency of the reduced collocation method is the ability to decompose the computation
into an offline component and an efficient online component. In this section, we describe how a complete
offline-online decomposition is achieved for the two algorithms. \rev{We also include an estimate of the computational
complexity,  which makes evident the dependence of the computational cost
on $\calN$ in the offline computation and its independence  in the online computation.}

\subsubsection{Least Squares}
\label{sec:OffOnLS}

We begin with the least-squares equation \eqref{eq:reducedprobLS},
\begin{equation*}
\mathbb{A}_N^T(\mu^*)\,\mathbb{A}_N(\mu^*) \,\vec{c} = \mathbb{A}_N^T(\mu^*)\,\boldf^\calN .
\end{equation*}
Invoking the affine assumption for  $\mathbb{L}$ (Equation \eqref{eq:affineL}) and $f$ (Equation \eqref{eq:affinef}) gives
\begin{eqnarray*}
\left(\mathbb{A}_N^T\,\mathbb{A}_N\right)_{ij} & =&
 \left(\mathbb{L}_\calN\, u_{\mu^i}^\calN\right)^T
\left(\mathbb{L}_\calN\, u_{\mu^j}^\calN\right) = \sum_{q_1 = 1}^{Q_a} \sum_{q_2 = 1}^{Q_a}
a_{q_1}^{\mathbb{L}}(\mu^*) a_{q_2}^{\mathbb{L}}(\mu^*) \left(\mathbb{L}_{q_1}\,
u_{\mu^i}^\calN\right)^T \left(\mathbb{L}_{q_2}\, u_{\mu^j}^\calN\right) \\
\left(\left(\mathbb{A}_N\right)^T\, \boldf^\calN \right)_i & = & \left(\mathbb{L}_\calN\, u_{\mu^i}^\calN\right)^T\, \boldf^\calN
 = \sum_{q_1 = 1}^{Q_a} \rev{\sum_{q_2 = 1}^{Q_f}} a_{q_1}^{\mathbb{L}}(\mu^*) \rev{a_{q_2}^f(\mu^*)}
 \left(\mathbb{L}_{q_1}\, u_{\mu^i}^\calN\right)^T \boldf_{q_2}^\calN.
\end{eqnarray*}
Hence, the decomposition \rev{and operation count} can be summarized as follows
\begin{itemize}
\item[{\bf Offline}] Calculate $\left(\mathbb{L}_{q_1}\,
u_{\mu^i}^\calN\right)^T \left(\mathbb{L}_{q_2}\, u_{\mu^j}^\calN\right)$ and $\left(\mathbb{L}_{q}\,
u_{\mu^i}^\calN\right)^T \boldf_{q_2}^\calN$ for $i,j \in \{1,\dots,N\}$, \rev{ with complexity of order
$N^2 Q_a^2 \calN^2 + N Q_a Q_f \calN$.}
\item[{\bf Online}] Form the $N \times N$ matrix $\mathbb{A}_N^T\,\mathbb{A}_N$ and $N \times 1$ vector
$\left(\mathbb{A}_N\right)^T\, \boldf^\calN$ \rev{for any $\mu^* \in \calD$}
 and solve the reduced $N \times N$ system for \rev{coefficients} $c_j$
 \eqref{eq:reducedprobLS}. \rev{Online complexity is of order $N^2 Q_a^2 + N^3 + Q_a Q_f N$.}
\end{itemize}

\subsubsection{Empirical Collocation}
\label{rmk:OffOnRC} Here, we demonstrate the offline-online decomposition for the reduced collocation
approach. The reduced equation in this case is
\begin{equation}
\sum_{j=1}^{N} c_j \mathbb{P}_\calN^N \left(\mathbb{L}_\calN(\mu^*) u^\calN_{\mu^j}\right)  = f(x; \mu^*).
\quad {\rm for} \quad x \in C_R^N
\end{equation}
which becomes \[
\sum_{j=1}^{N} c_j \sum_{q = 1}^{Q_a} a^{\mathbb{L}}_q(\mu^*) \mathbb{P}_\calN^N \left(\mathbb{L}_q
u^\calN_{\mu^j}\right) = \sum_{q = 1}^{Q_f} a^f_q(\mu^*) f_q(x)
\]
with the affine assumptions \eqref{eq:affineL} and \eqref{eq:affinef}.

This means that, given $\{\mu^1, \mu^2, \dots, \mu^N\}$ and the set of $N$ reduced collocation points
$C^N_R$, the splitting of the computation is done as follows:
\begin{itemize}
\item[{\bf Offline}] Calculate $\mathbb{L}_q u^\calN_{\mu^j}$, their $\mathbb{P}_\calN^N$ projections, and
$f_q(x)$ \rev{for $x \in C^N_R$}. \rev{The complexity is of order $N^2 Q_a \calN^3$ (see Section
\ref{sec:Proj} for the complexity for the projection $\mathbb{P}_\calN^N$).}
\item[{\bf Online}] Form $a^{\mathbb{L}}_q(\mu^*) \mathbb{P}_\calN^N \left(\mathbb{L}_q u^\calN_{\mu^j}\right)$ for any $j$ and $q$,
evaluate $a^f_q(\mu^*)$ and form $a^f_q(\mu^*) f_q(x)$ at the reduced set of collocation points $C_R^N$, and
finally solve the reduced $N \times N$ system for $c_j$'s \eqref{eq:rbmpderc}. \rev{Online complexity is
of the order $Q_a N^2 + N^3 + Q_f$.}
\end{itemize}

\subsection{Efficiently computing the Error Estimator}
\label{sec:OffOnEE} Although we are primarily interested in minimizing the online cost of computation, it is
also advantageous to be able to efficiently compute the offline component of the reduced collocation method. In
particular, the  greedy algorithm requires repeated computations of the error estimator $\Delta_i(\mu)$ for
$i \in \{1, \dots, N\}$ and any $\mu \in \Xi$.
 To make this practical, as we select more and more bases and $i$  goes from $1$ to $N$, we can reuse previously computed components of the error estimator.
This can be achieved in essentially the same fashion as in the Galerkin framework. Indeed, we have,
\begin{alignat*}{1}
\lVert \rev{\boldf^\calN - \mathbb{L}_\calN (\mu) u_\mu^{(i)}}\rVert_{\ell^2}^2 & = (\boldf^\calN-
\mathbb{L}_\calN (\mu) u_\mu^{(i)})^T (\boldf^\calN - \mathbb{L}_\calN (\mu) u_\mu^{(i)}).
\end{alignat*}
%\mathbb{L}_\calN (\mu)^T \mathbb{L}_\calN (\mu)
The resulting three terms after expansion are
\begin{alignat*}{1}
e^i_1(\mu) & := (\boldf^\calN)^T  \boldf^\calN,\\
e^i_2(\mu) & :=(\mathbb{L}_\calN (\mu) u_\mu^{(i)})^T (\mathbb{L}_\calN (\mu)
u_\mu^{(i)}),\\
e^i_3(\mu) & :=(\boldf^\calN)^T (\mathbb{L}_\calN (\mu) u_\mu^{(i)}).
\end{alignat*}
They can be handled efficiently in the same fashion. To do that, we invoke the affine assumptions
\eqref{eq:affineL}-\eqref{eq:affinef} and the expansion of the reduced solution $u_\mu^{(i)} =
\sum_{j=1}^i\,c_j(\mu) u_{\mu^j}^\calN$ to obtain
\begin{alignat*}{1}
e^i_1(\mu) & = \sum_{q_3,q_4=1}^{Q_f}
a_{q_3}^f(\mu)a_{q_4}^f(\mu)\,f_{q_3}^T f_{q_4}\\
e^i_2(\mu) & = \sum_{j_1,j_2=1}^{i}\sum_{q_1,q_2=1}^{Q_a}c_{j_1}(\mu)c_{j_2}(\mu)
a_{q_1}^{\mathbb{L}}(\mu)a_{q_2}^{\mathbb{L}}(\mu)\,\left(u_{j_1}^\calN\right)^T \mathbb{L}_{q_1}^T \mathbb{L}_{q_2} u_{j_2}^\calN\\
e^i_3(\mu) & = \sum_{q_1=1}^{Q_a} \sum_{q_3=1}^{Q_f}\sum_{j_1=1}^{i} a_{q_1}^{\mathbb{L}}(\mu)
a_{q_3}^f(\mu) c_{j_1}(\mu)\,f_{q_3}^T \mathbb{L}_{q_1} u_{j_1}^\calN
\end{alignat*}
The Offline-Online decomposition of these terms \rev{and their computational complexities are} as follows.
\begin{itemize}
\item[{\bf Offline}] Calculate
\begin{alignat*}{1}
&f_{q_3}^T f_{q_4},\\
&\left(u_{j_1}^\calN\right)^T \mathbb{L}_{q_1}^T \mathbb{L}_{q_2}u_{j_2}^\calN,\\
&f_{q_3}^T \mathbb{L}_{q_1} u_{j_1}^\calN
\end{alignat*}
for $q_1,q_2 \in \{1,\dots, Q_a\}$, \,\, $q_3,q_4 \in \{1,\dots, Q_f\}$ and $j_1,j_2 \in \{1, \dots, i\}$.
\rev{The cost is of order $Q_f^2 \calN + N^2 \calN^2 + Q_a Q_f N \calN^2$.}
\item[{\bf Online}] Evaluate the coefficients
\begin{alignat*}{1}
&a_{q_3}^f(\mu)a_{q_4}^f(\mu),\\
&c_{j_1}(\mu)c_{j_2}(\mu)
a_{q_1}^{\mathbb{L}}(\mu)a_{q_2}^{\mathbb{L}}(\mu),\\
&a_{q_1}^{\mathbb{L}}(\mu) a_{q_3}^f(\mu) c_{j_1}(\mu),
\end{alignat*}
and form
\begin{alignat*}{1}
e^i_1(\mu) & := (\boldf^\calN)^T  \boldf^\calN,\\
e^i_2(\mu) & :=(\mathbb{L}_\calN (\mu) u_\mu^{(i)})^T (\mathbb{L}_\calN (\mu)
u_\mu^{(i)}),\\
e^i_3(\mu) & :=(\boldf^\calN)^T (\mathbb{L}_\calN (\mu) u_\mu^{(i)}).
\end{alignat*}
\rev{The online computation has complexity of order $Q_f^2 + N^2 Q_a^2 + Q_a Q_f N$.}
%$e^i_1(\mu) = f^T \mathbb{L}_\calN (\mu)^T \mathbb{L}_\calN (\mu)f$, \,\, $e^i_2(\mu) =(\mathbb{L}_\calN
%(\mu) u_\mu^N)^T \mathbb{L}_\calN (\mu)^T \mathbb{L}_\calN (\mu) (\mathbb{L}_\calN (\mu) u_\mu^N)$, \,\,
%$e^i_3(\mu) =f^T \mathbb{L}_\calN (\mu)^T \mathbb{L}_\calN (\mu) (\mathbb{L}_\calN (\mu) u_\mu^N)$.
\end{itemize}

\subsection{Comparison with the Galerkin RBM}
\label{sec:comparison_rc}

In this section, we show a particular advantage of the proposed Empirical Reduced Collocation Method over
the traditional reduced basis method in the Galerkin framework. When the operator is non-affine, that is, we
have instead of \eqref{eq:affineL}
\begin{equation}
\mathbb{L} (\mu) = \sum_{q = 1}^{Q_a} a^{\mathbb{L}}_q(x,\mu) \mathbb{L}_q, \label{eq:nonaffineL}
\end{equation}
The Galerkin approach has to use the Empirical Interpolation Method
\cite{Barrault_Nguyen_Maday_Patera,Grepl_Maday_Nguyen_Patera} to achieve the offline-online decomposition
and the traditional speedup. In fact, $a^{\mathbb{L}}_q(x,\mu)$ has to be approximated by the affine
expansion
\begin{equation}
\label{eq:nonaffine_expL} a^{\mathbb{L}}_q(x,\mu) = \sum_{m=1}^{M_q}\,\phi^q_m (\mu) a_{aff,m}^q (x),
\end{equation}
so that $(\mathbb{L}_q u^\calN_{\mu^i},u^\calN_{\mu^j}a_{aff,m}^q (x))_\Omega$ are computed offline for all
$i,j,q,m$. During the online stage for any given $\mu$, $\phi^q_m (\mu)$ are obtained and
$\sum_{q=1}^{Q_a}\sum_{m=1}^{M_q}(\mathbb{L}_q u^\calN_{\mu^i},u^\calN_{\mu^j}a_{aff,m}^q (x))_\Omega$ are
formed. Obviously, the online performance is dependent on
$ \sum_{q=1}^{Q_a} M_q.$
The proliferation from $Q_a$ to $\sum_{q=1}^{Q_a} M_q$ adversely affects the online performance of the reduced
basis method and limits its practical scope. This is particularly the case for geometrically complex
problems with parameters describing the geometry \cite{CHM_JCP, PomplunSchmidt2010, Rozza_Huynh_Patera}:
$\sum_{q=1}^{Q_a} M_q$ can be one to two magnitudes larger than $Q_a$. The online efficiency is thus
significantly worse than the affine problems.

However, this significant barrier {\em does not exist} for the proposed empirical reduced collocation
method. Since to form the online solver we only need to evaluate $a^{\mathbb{L}}_q(x,\mu^*)$ for $x \in
C_R^N$. This can be done without the expansion \eqref{eq:nonaffine_expL}. Note that $\mathbb{P}_\calN^N
\left(\mathbb{L}_q u^\calN_{\mu^j}\right)$ is readily available from the offline calculation.

Unfortunately, this advantage of the empirical reduced collocation method over the Galerkin framework does
not translate to the least squares reduced collocation method: when $a_q(\mu) = a_q(x,\mu)$, we need to
perform the expansion \eqref{eq:nonaffine_expL} to have the online procedure of forming $N \times N$ matrix
$\mathbb{A}_N^T\,\mathbb{A}_N$ independent of $\calN$. The fundamental reason is that least squares is
intrinsically a projection method and thus our least squares reduced collocation method is closely related to
the Galerkin RBM framework.

\section{Numerical Results}
\label{sec:numerical}

In this section, we consider a couple of two-dimensional diffusion-type test problems similar to those used in
\cite{Prudhomme_Rovas_Veroy_Maday_Patera_Turinici,Rozza_Huynh_Patera} to show the accuracy and efficiency of
the proposed methods:
\begin{enumerate}
\item Diffusion
\begin{eqnarray}
(1+ \mu_1 x) u_{xx} + (1+ \mu_2 y) u_{yy}= e^{4 x y}
\end{eqnarray}
on $\Omega = [-1,1] \times [-1,1]$ with zero Dirichlet boundary condition.
%\item Convection Diffusion
\item Anisotropic wavespeed simulation
\begin{eqnarray}
-u_{xx} - \mu_1 u_{yy} -\mu_2 u = -10\,\sin(8\,x\,(y-1))
\end{eqnarray}
on $\Omega = [-1,1] \times [-1,1]$ with zero Dirichlet boundary condition.
\end{enumerate}

Our truth approximations are generated by a spectral Chebyshev collocation method
\cite{TrefethenSpecBook,HesthavenGottlieb2007}. For $C^\calN$ and for $X$, we use the Chebyshev  grid based
on \rev{$\calN_x$} points in each direction \rev{with $\calN_x^n = \calN$}. We consider the parameter domain $\calD$ for $(\mu_1,\mu_2)$  to be
$[-0.99,0.99]^2$ and $[0.1,4] \times [0,2]$ respectively for the two test problems. For $\Xi$, they are
discretized uniformly by $64 \times 64$ and $128 \times 64$ Cartesian grids. \rev{For the purpose of reproducible
research, the code has been posted online \cite{YChenSite}.}

\subsection{Preliminaries}

\subsubsection{Computation of $\mathbb{P}_\calN^N$}
\label{sec:Proj}

For the empirical reduced collocation method, we need the fine-to-coarse projection $\mathbb{P}_\calN^N$. We
begin with a set of Chebyshev points in one dimension $x_j = \cos(\frac{\pi j}{\calN})$. Given a vector of
function values $w(x_j)$, we define the function $w(x)$ by the Chebyshev expansion
\cite{HesthavenGottlieb2007}
\begin{equation}
w(x) = \sum_{k=0}^\calN a_k T_k(x). \label{eq:chebexp}
\end{equation}
where
\begin{equation} a_k =  \frac{2}{\calN c_k}\sum_{j = 0}^\calN\frac{1}{c_j} w(x_j)
\cos(k\,\arccos(x_j)) = \frac{2}{\calN c_k}\sum_{j = 0}^\calN\frac{1}{c_j} w(x_j) \cos \left( \frac{\pi j
k}{\calN} \right).\label{eq:chebycoe}\end{equation} Here,
\[
c_k =
\begin{cases}
2 \quad & k = 0, \calN\\
1 \quad & \mbox{ otherwise }
\end{cases}.
\]

This definition relies on the fact that  the Chebyshev polynomial is
\[ T_k(x) = \cos( k \arccos(x)),\]
so that
 \[w(x_j)= \sum_{k=0}^\calN a_k  \cos(\frac{\pi j k}{\calN}) . \]

Now, if we wish to evaluate the function value of $w(x)$ at any set of points $\{x_\ell \}$, we simply  plug
those points into the Chebyshev expansion
\begin{equation}
w(x_\ell) = \sum_{k=0}^\calN a_k T_k(x_\ell). \label{eq:chebeval}
\end{equation}
 In particular, the calculation of $\mathbb{P}_\calN^N w$ is done by evaluation
at the reduced set of $N$ points $C_R^N$.
\rev{ For two or three space dimensions,   \eqref{eq:chebycoe} becomes a double or  triple sum
over the points in each dimension,  while \eqref{eq:chebeval} contains only the reduced collocation
points. Thus,  computational complexity to
evaluate \eqref{eq:chebeval} for $\ell \in \{1, \dots, N\}$ is of order $N \calN^2$, where $\calN$ is the total number
of points in the fine mesh and $N$ is the number of points in the reduced collocation grid.
Unfortunately,  $\calN$ is a product of  the number of points in each dimension, which of course grows
exponentially.}
% Now if we wish to find the values $w^N(y_j)$, where $y_j =  \cos(\frac{\pi j}{M_{rb}})$, it is also quite simple because
%  \[ v_j = w^N(y_j)= \sum_{k=0}^N a_k  \cos(\frac{\pi j k}{M_{rb}}) , \; \; \; \; j=0, . . . M_{rb} .\]
%This vector has length $M_{rb}$, and now we solve the system
%\[ \sum_{j=1}^{M_{rb}} c_j v_j = e^{4 y} \]

\subsubsection{Results of the fine solver and setup for the reduced solvers}

\begin{figure}[ht]
  \begin{center}
        \includegraphics[width=\textwidth]{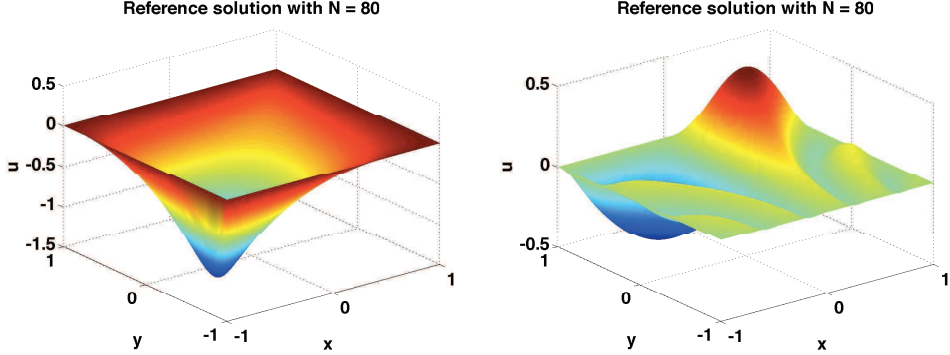}
  \end{center}
    \caption{The truth approximations for diffusion (left) and anisotropic wavespeed simulation (right) for $\mu_1 = 1$ and $\mu_2 = 0.5$ computed on a $81 \times 81$ Chebyshev grid.}
    \label{fig:truthsample}
\end{figure}

\begin{figure}[ht]
  \begin{center}
        \includegraphics[width=0.49\textwidth]{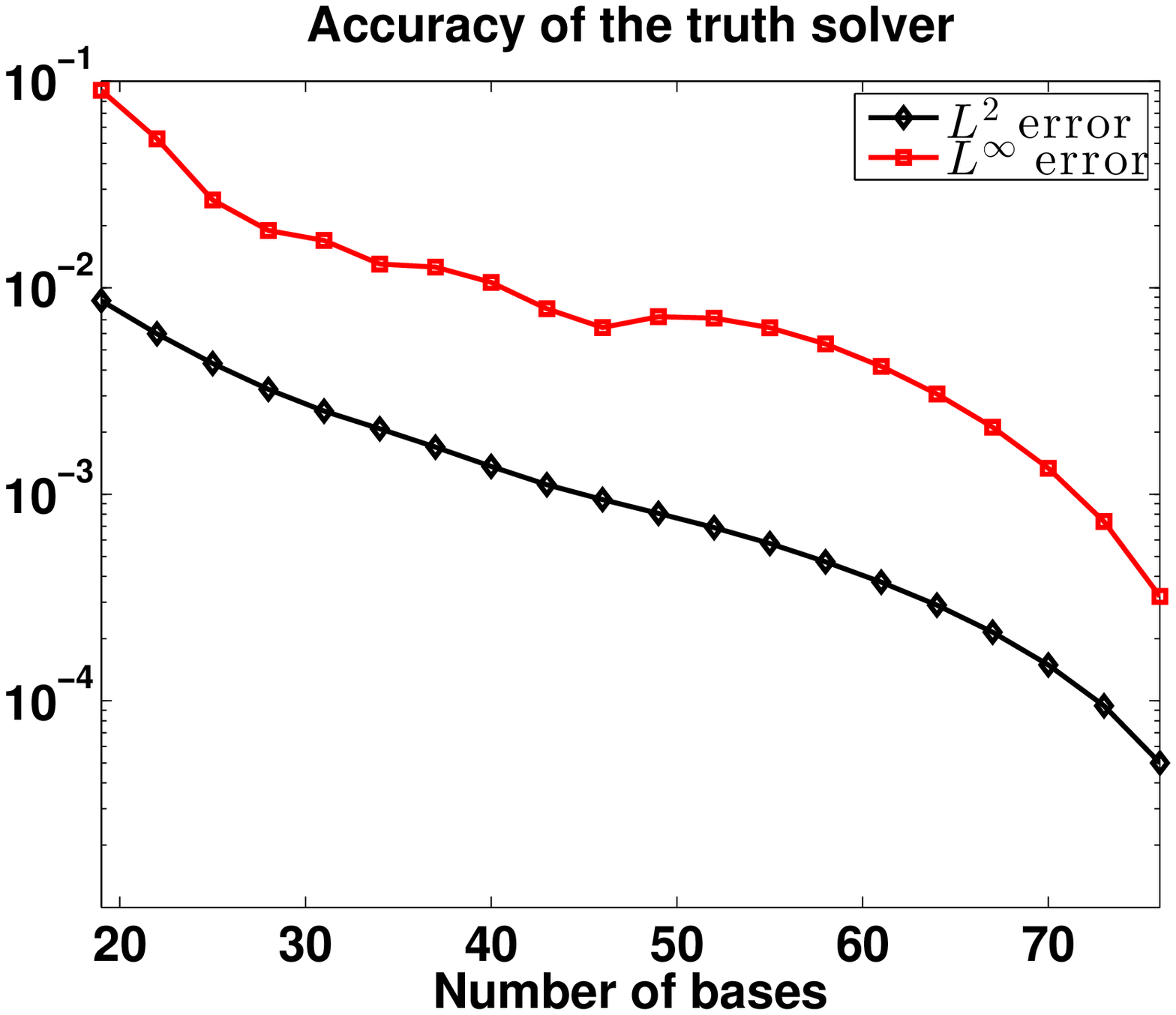}
        \includegraphics[width=0.49\textwidth]{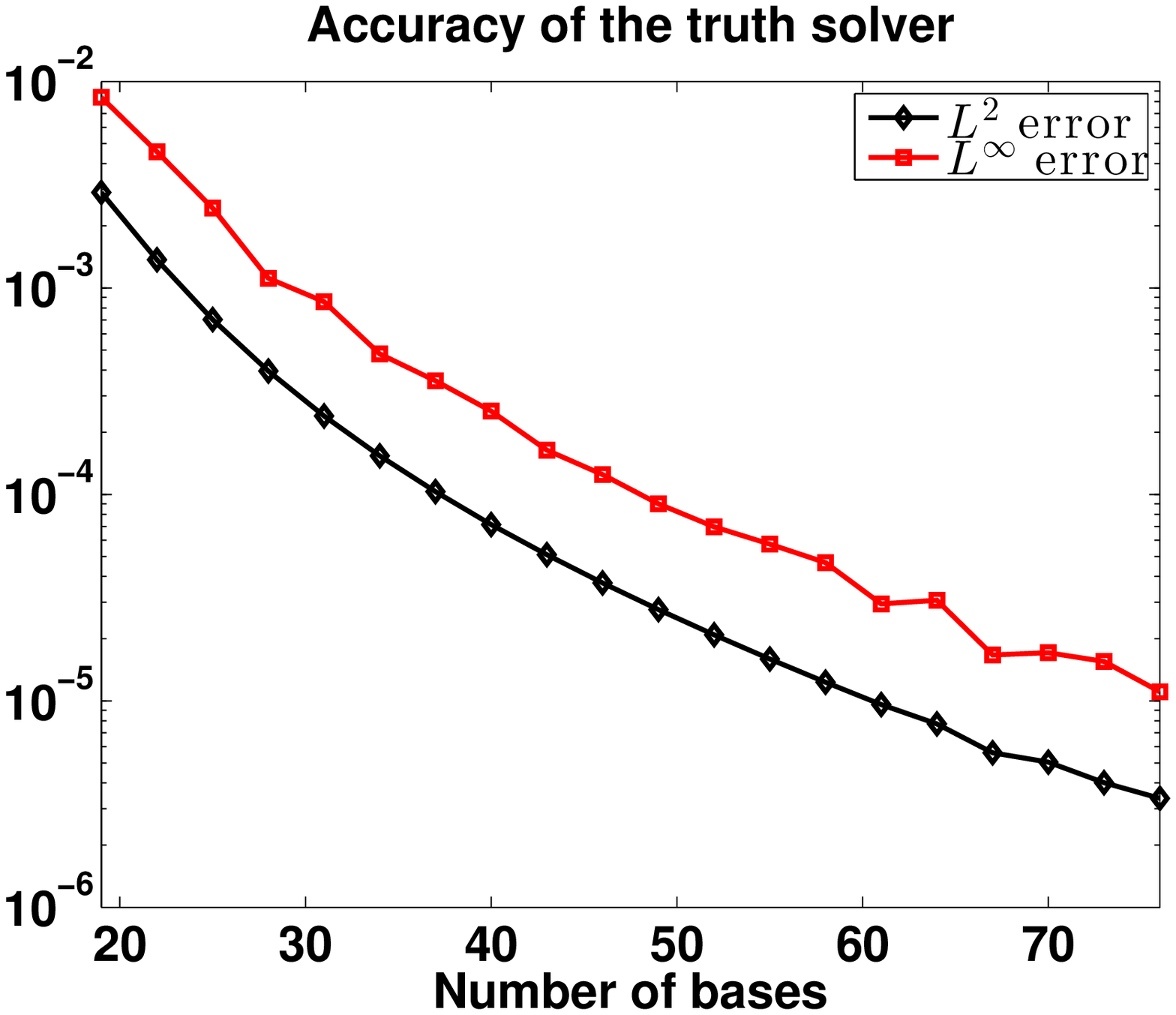}
  \end{center}
    \caption{The accuracy of truth approximations for diffusion (left) and anisotropic wavespeed
    simulation (right) for $\mu_1 = 1$ and $\mu_2 = 0.5$ with reference solution being computed on
    a $81 \times 81$ Chebyshev grid.}
    \label{fig:truthaccuracy}
\end{figure}

\begin{figure}[ht]
  \begin{center}
        \includegraphics[width=\textwidth]{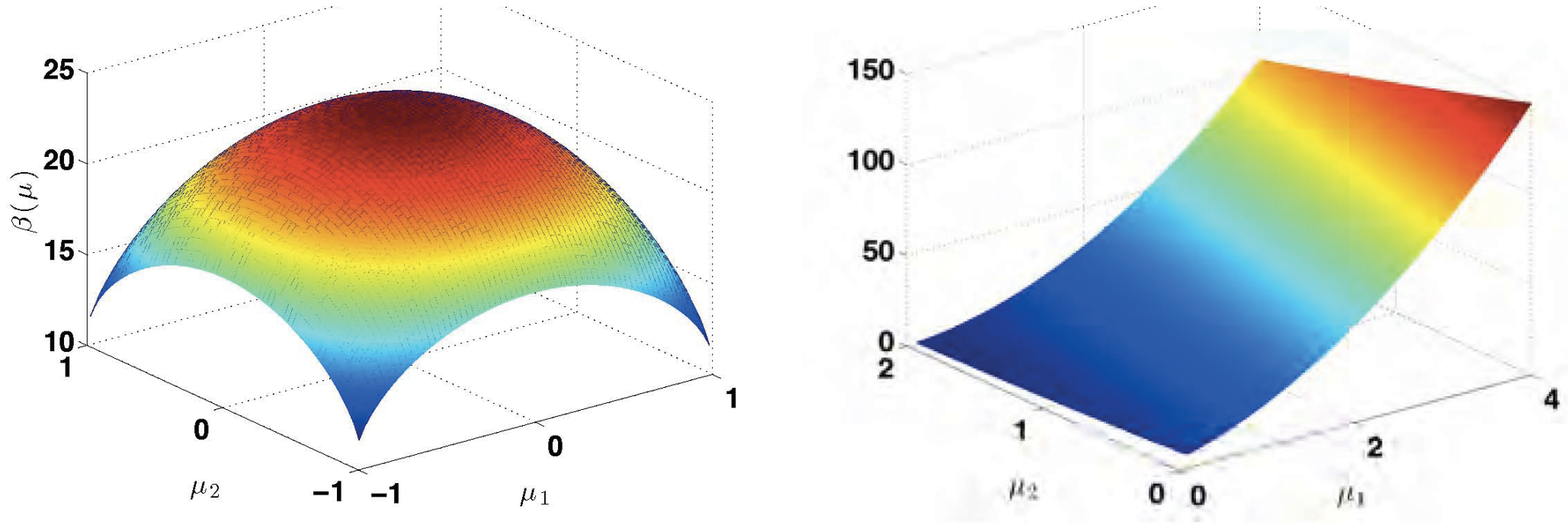}
  \end{center}
    \caption{The eigenvalues for the diffusion problem (left) and Anisotropic Wavespeed Simulation (right).}
    \label{fig:diffanisoeigs}
\end{figure}

Before we begin with the reduced basis solver, we must quantify the accuracy of the fine-domain solver,
which produces the truth approximations. Reference solutions computed by Chebyshev
collocation method on a $81 \times 81$ grid for $\mu_1 = 1$ and $\mu_2 = 0.5$ are plotted in Figure
\ref{fig:truthsample}. We also compute the truth solutions on a $\calN_x \times \calN_x$ grid for $\calN_x$
changing from $20$ to $77$ and evaluate the $L^2$ and $L^\infty$ errors. Exponential convergence of the
truth approximation with respect to $\calN$ is shown by Figure \ref{fig:truthaccuracy} as expected.

In the greedy algorithm, we required a lower bound on the eigenvalue of the operator. For the purposes of
this work, we simply calculate the smallest eigenvalue for each $\mu \in \Xi$ and use it as the lower bound
$\beta_{LB}(\mu)$. \rev{That is, the error estimator for the reduced basis solution $u_\mu^{(i)}$ based on
a reduced basis set $\{u^\calN_{\mu^1},\dots,u^\calN_{\mu^i}\}$ is given by
\[
\Delta_i(\mu) = \frac{\lVert {\boldf^\calN - \mathbb{L}_\calN (\mu)
u_\mu^{(i)}}\rVert_{\ell^2}}{{\sqrt{\beta(\mu)}}}.
\]
}There are more efficient ways \cite{HuynhSCM,CHMR-Cras,HKCHP}. However, algorithm design and implementation
of how to efficiently calculate $\beta_{LB}(\mu)$ is not an emphasis of this paper. Instead, we are
concentrating on the design of the overall reduced basis method in the collocation framework.  The
eigenvalues $\beta(\mu)$ are plotted in Figure \ref{fig:diffanisoeigs} for the two test problems. The first
problem becomes close to being degenerate at the four corners of the parameter domain.

\subsection{Results of the reduced solver: Anisotropic wavespeed simulation}

In this section, we present the results of the two reduced collocation methods applied to the anisotropic
wavespeed simulation.

\begin{figure}[ht]
  \begin{center}
        \includegraphics[width=0.49\textwidth]{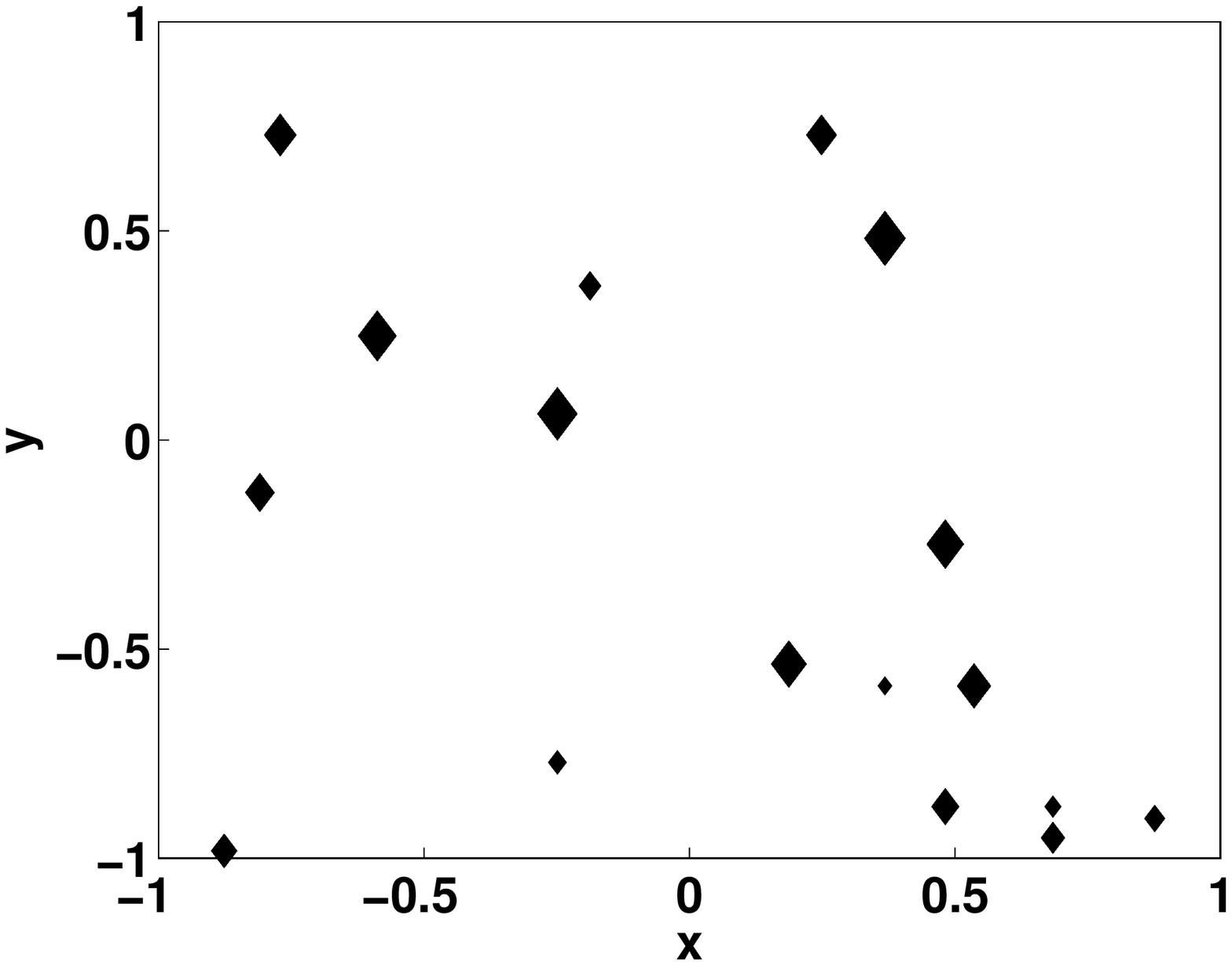}
        \includegraphics[width=0.49\textwidth]{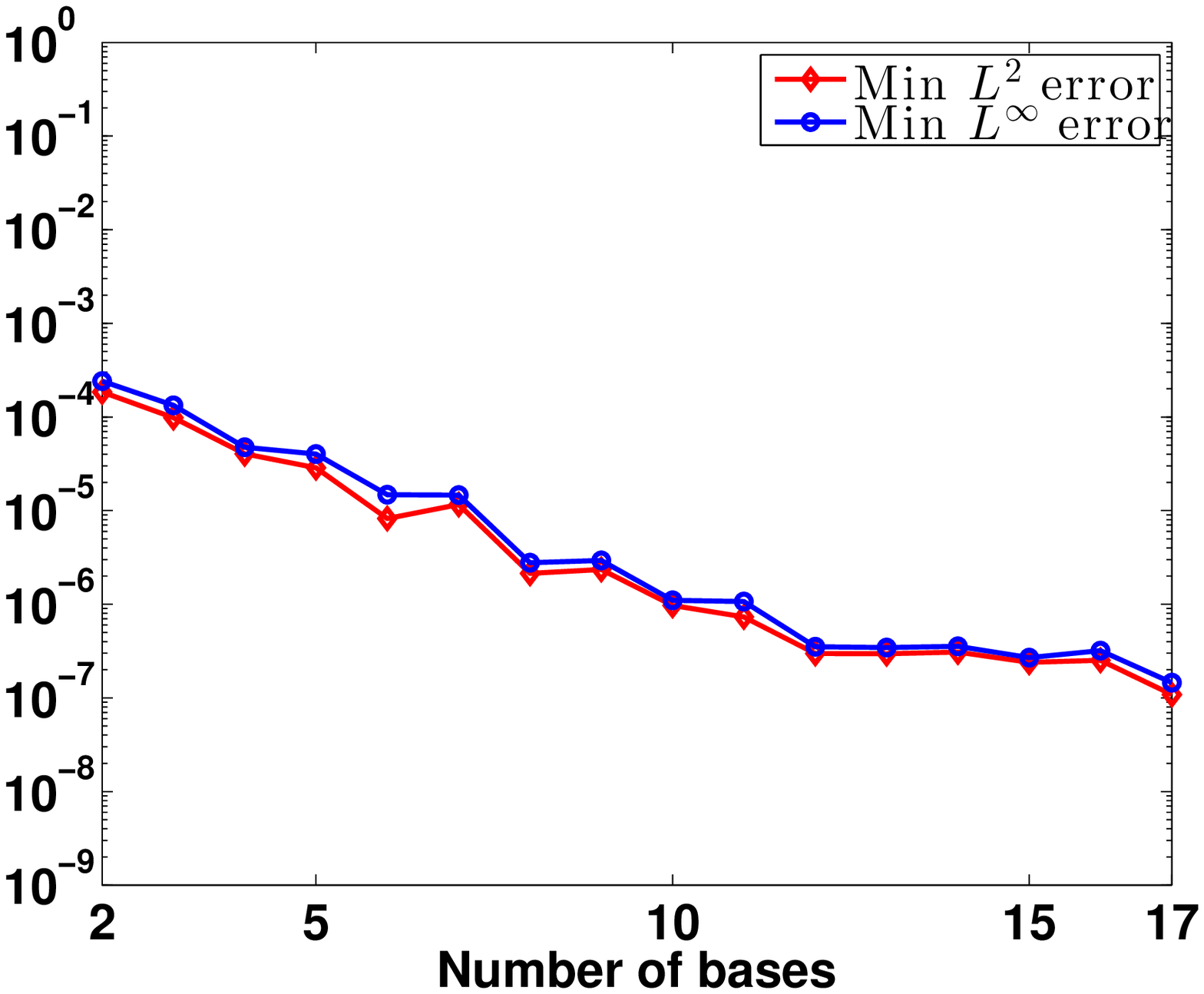}
        \includegraphics[width=0.49\textwidth]{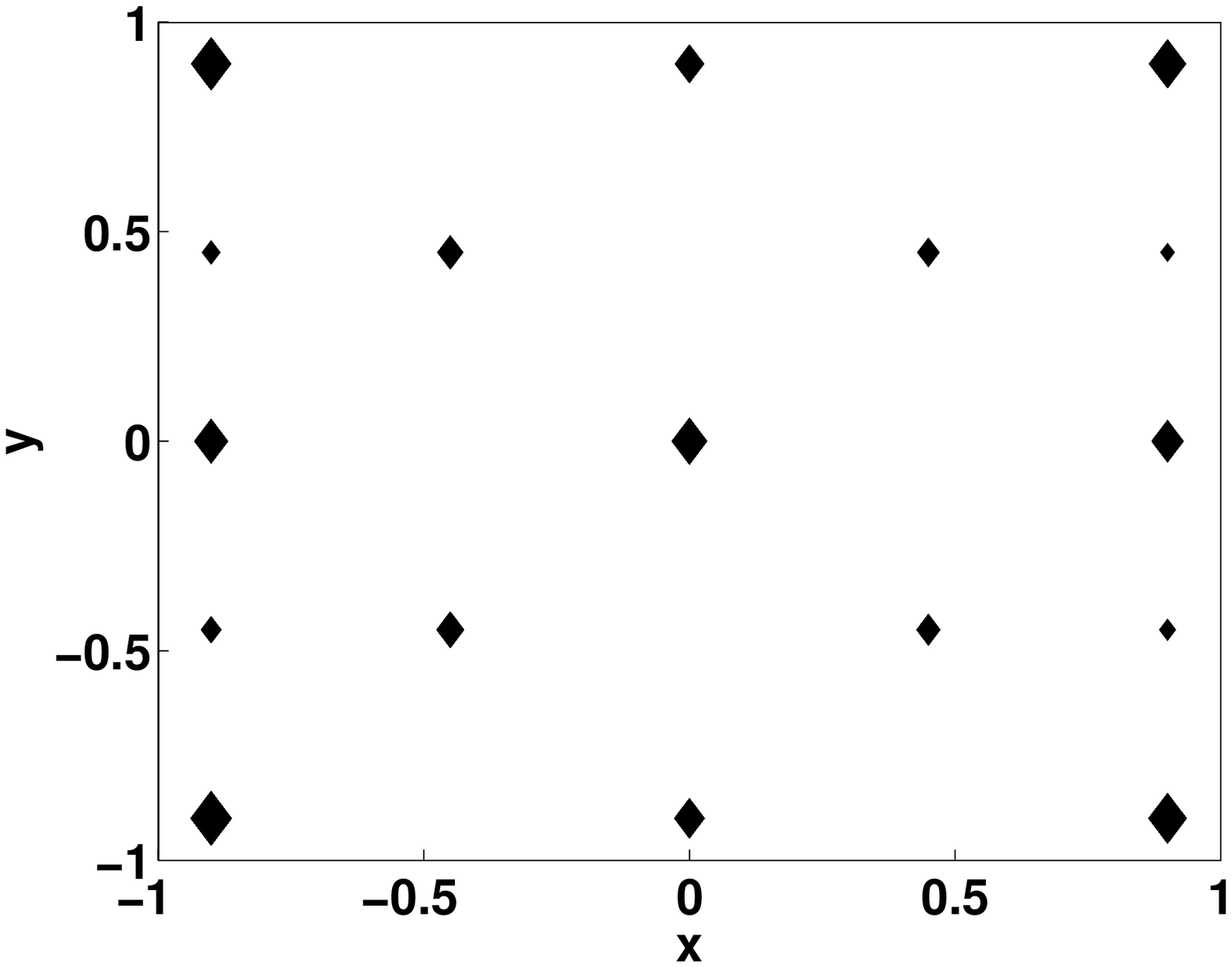}
        \includegraphics[width=0.49\textwidth]{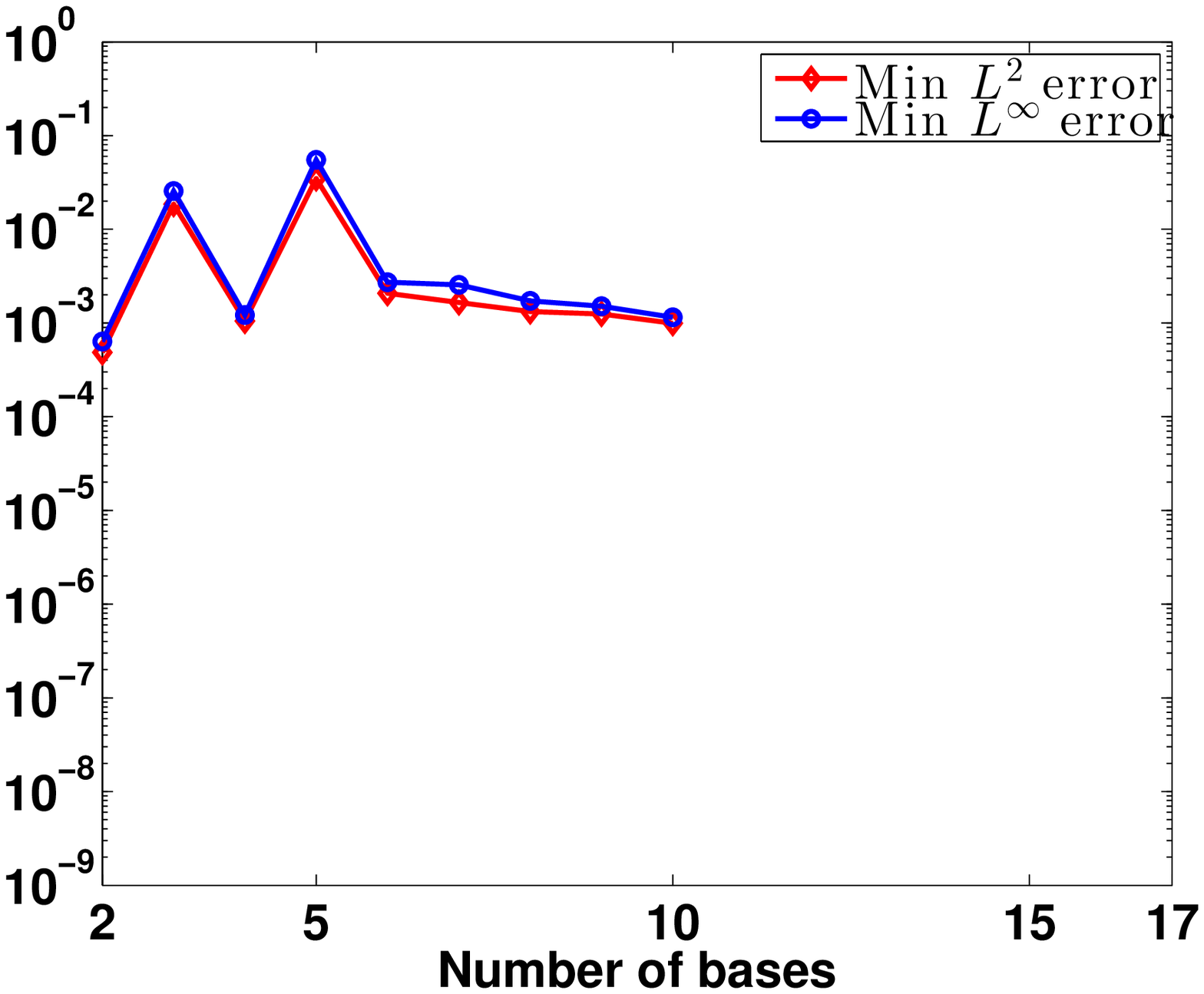}
  \end{center}
    \caption{Top: The reduced collocation points selected by the ERCM for Anisotropic wavespeed simulation, and
    the corresponding convergence plot. Bottom: The result of a (pre-determined) coarse Chebyshev grid was used.}
    \label{fig:diffecmAniso}
\end{figure}

\begin{figure}[ht]
  \begin{center}
        \includegraphics[width=\textwidth]{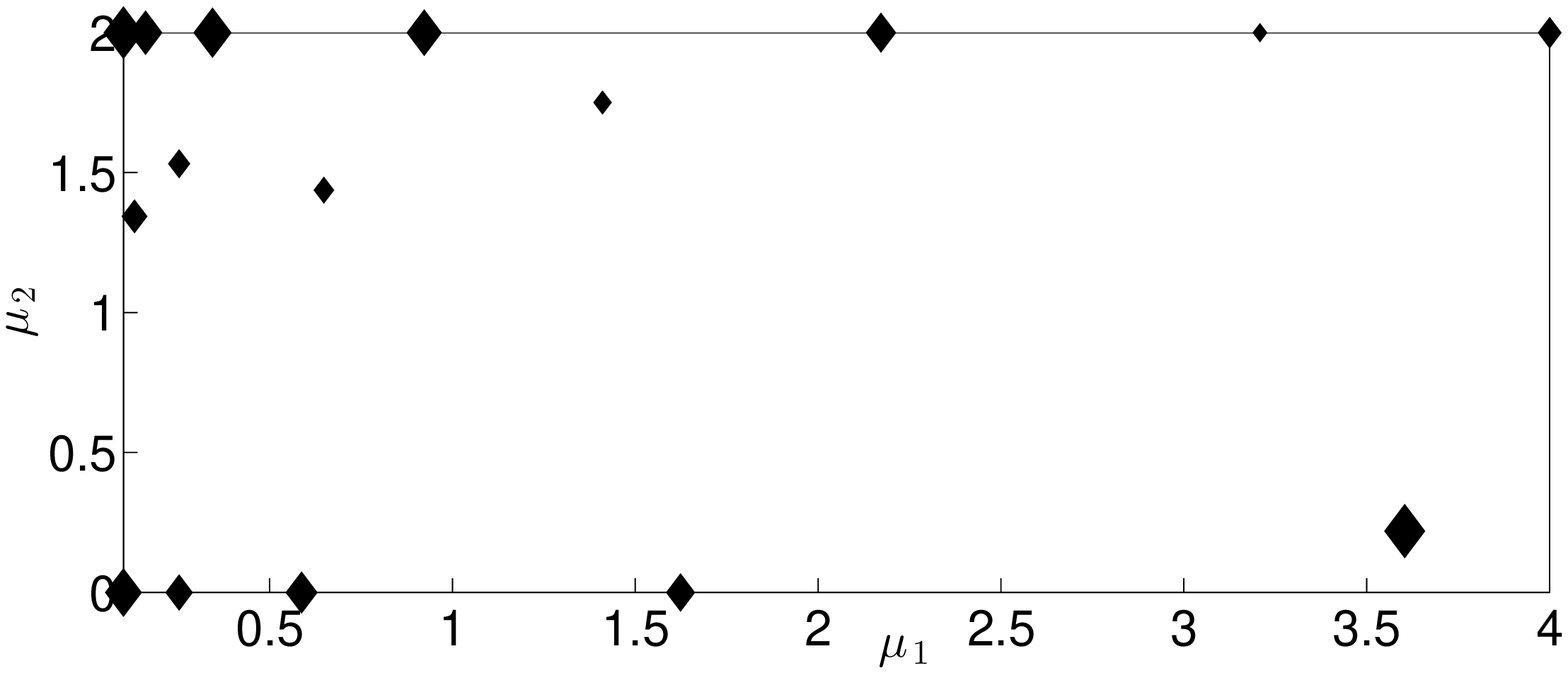}
        \includegraphics[width=\textwidth]{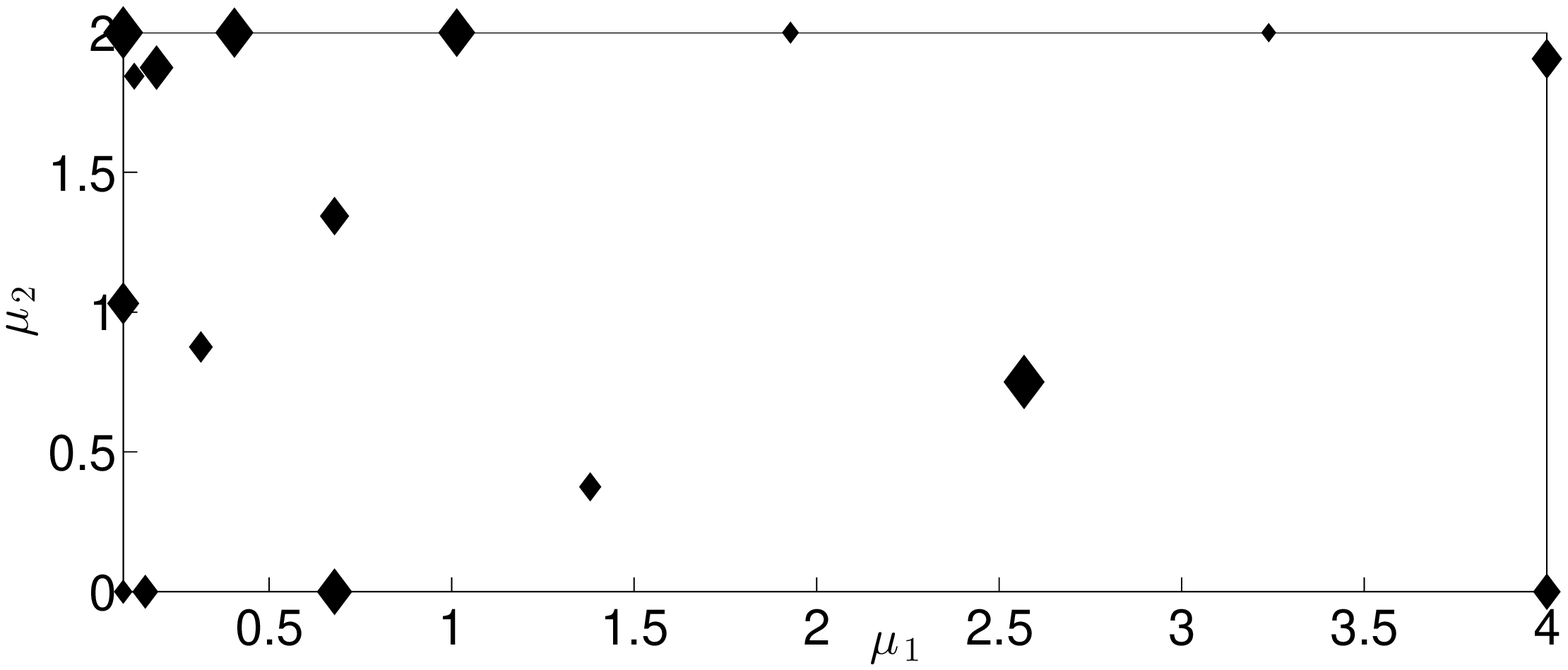}
  \end{center}
    \caption{The parameters picked by the greedy algorithm for pre-computation for the anisotropic wavespeed simulation. Top: Least Squares approach,
    bottom: Empirical Reduced Collocation Method.}
    \label{fig:pickedmu}
\end{figure}

\begin{figure}[ht]
  \begin{center}
        \includegraphics[width=0.45\textwidth]{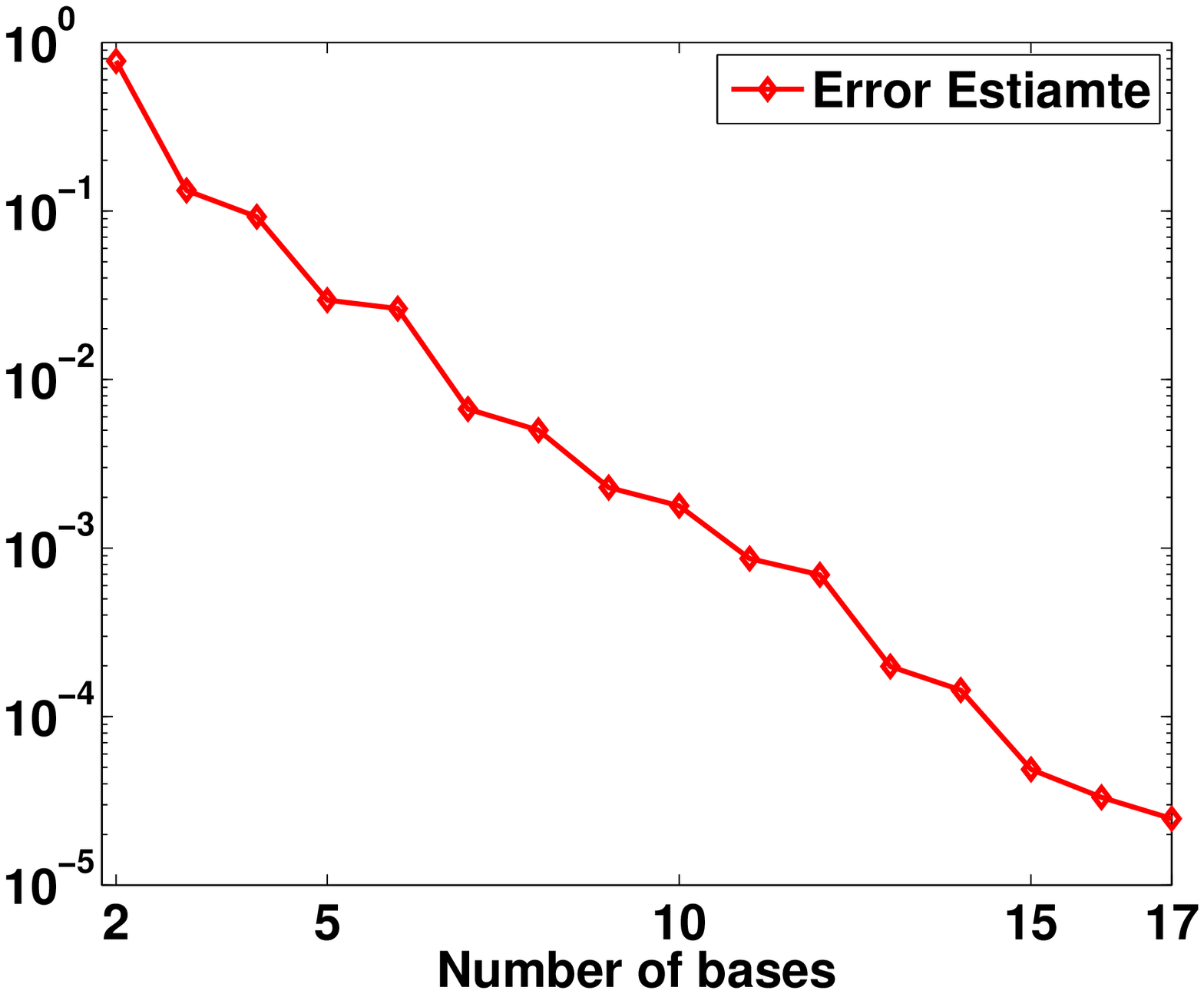}
        \includegraphics[width=0.45\textwidth]{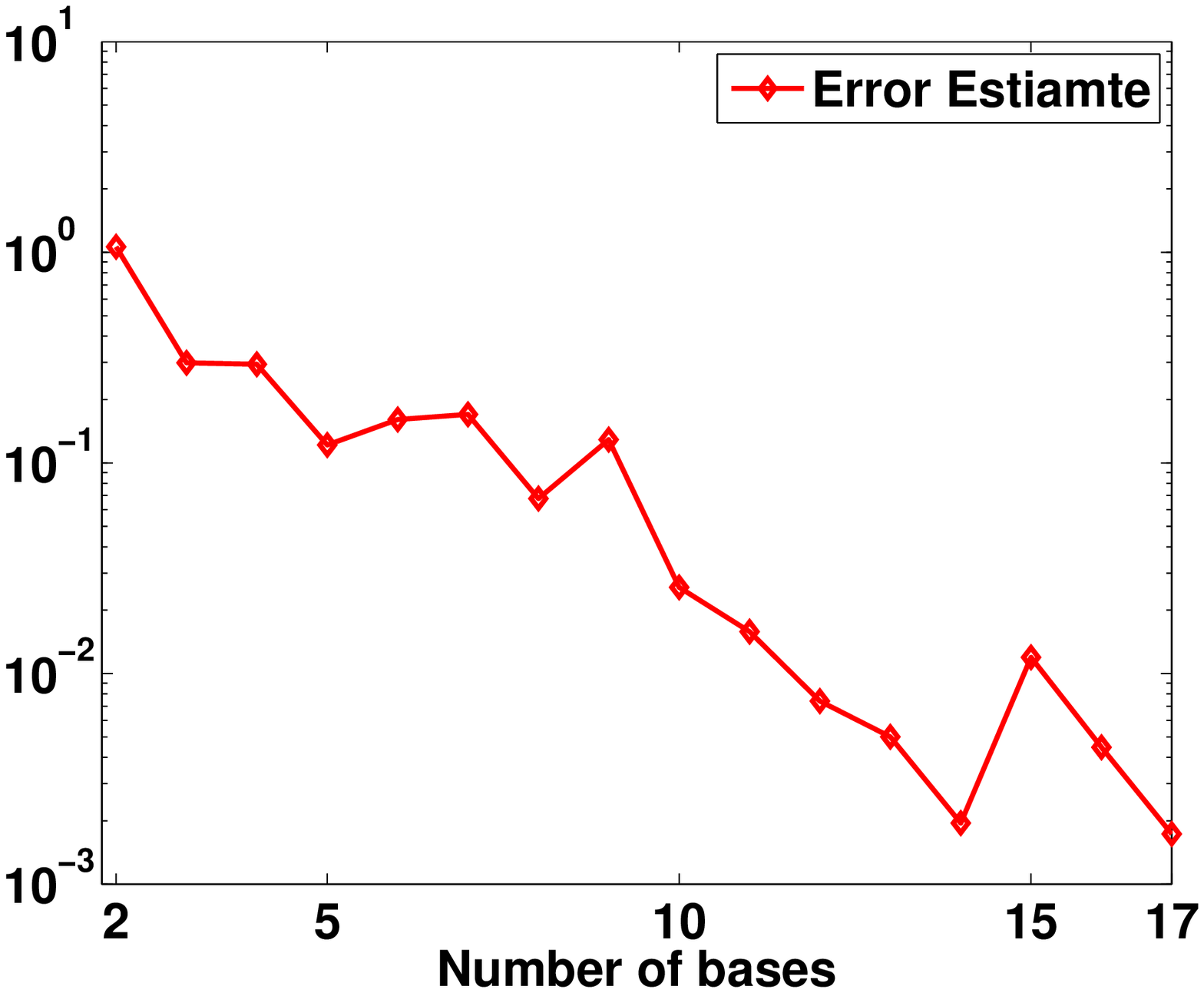}\\
        \includegraphics[width=0.45\textwidth]{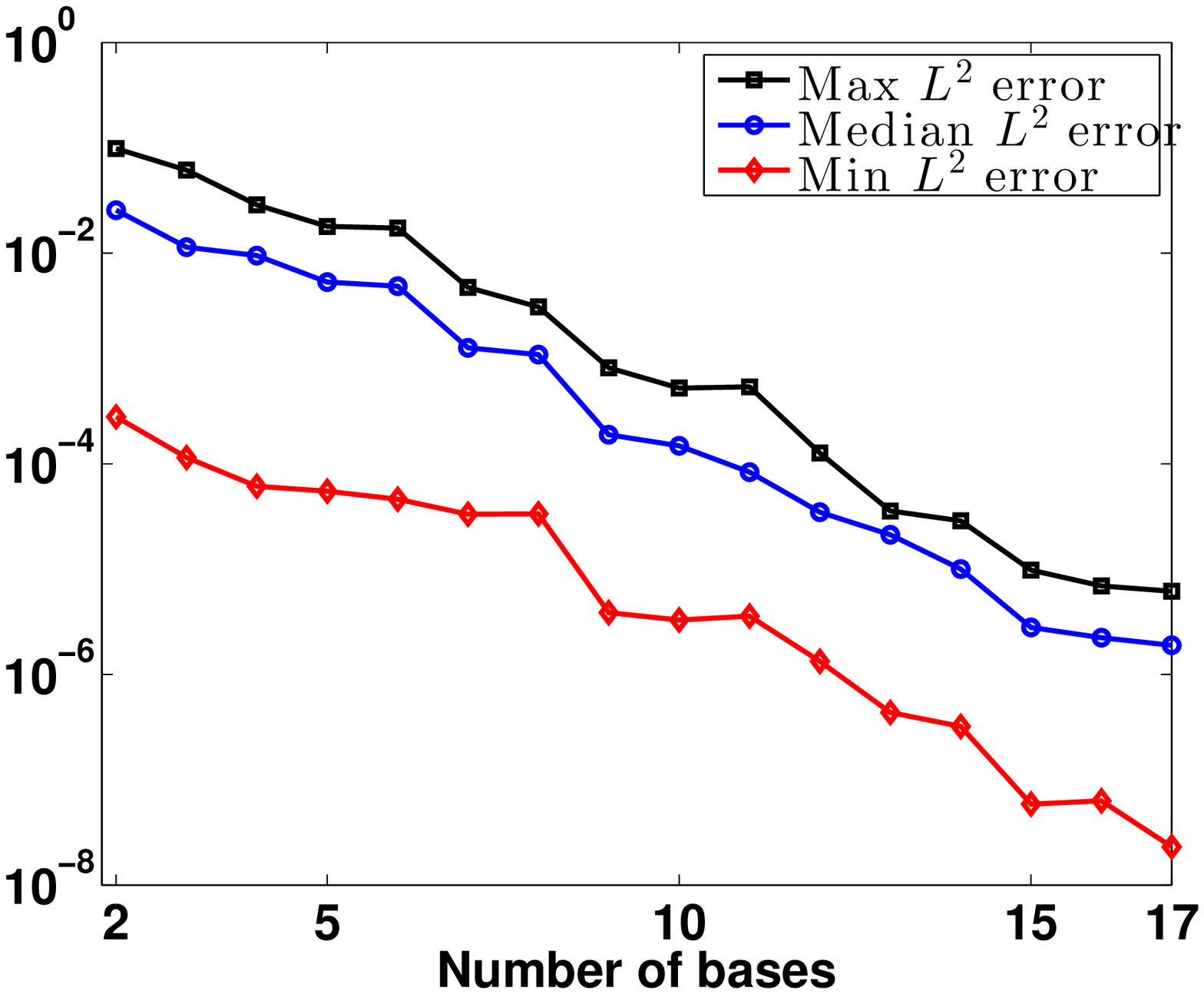}
        \includegraphics[width=0.45\textwidth]{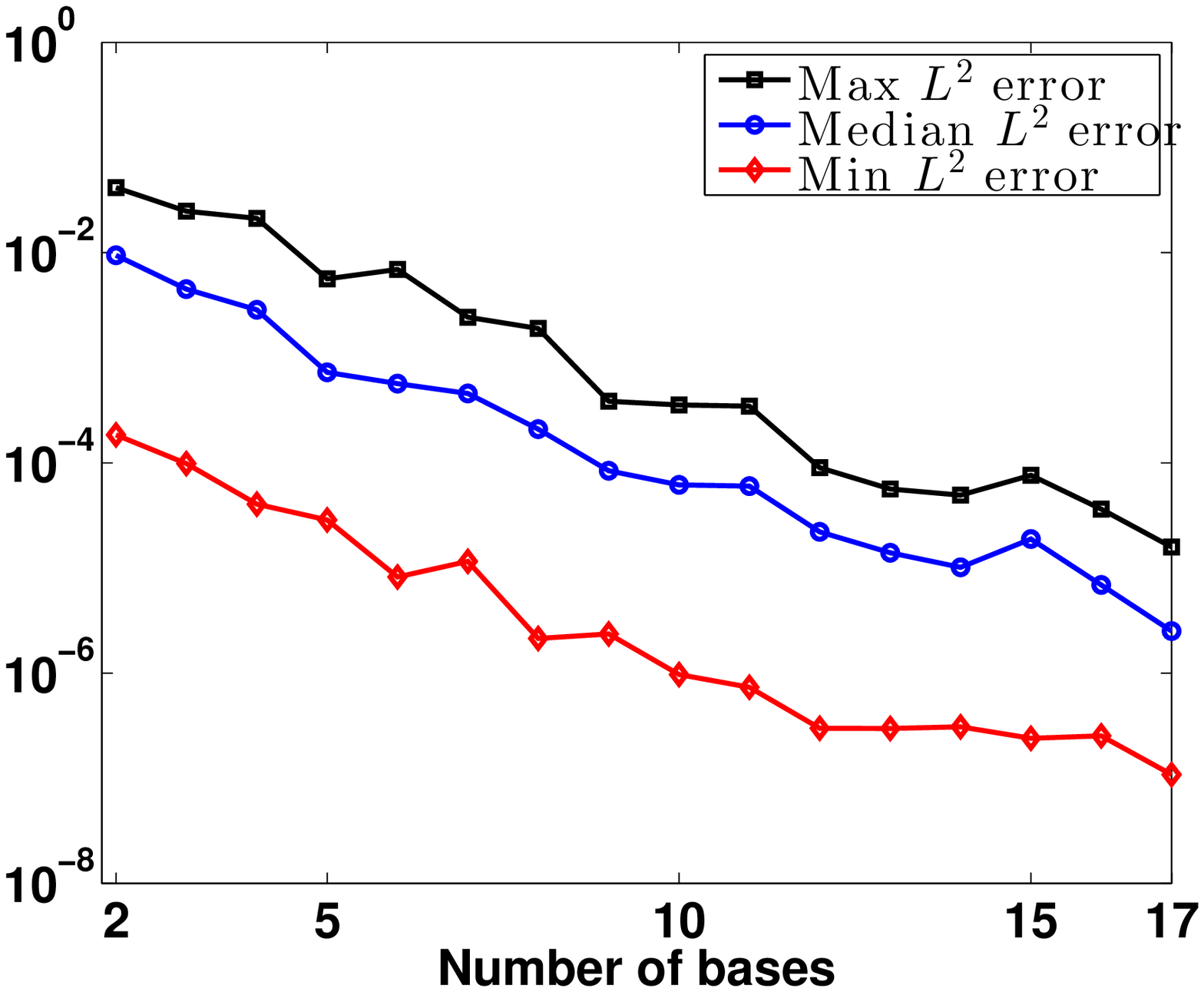}\\
        \includegraphics[width=0.45\textwidth]{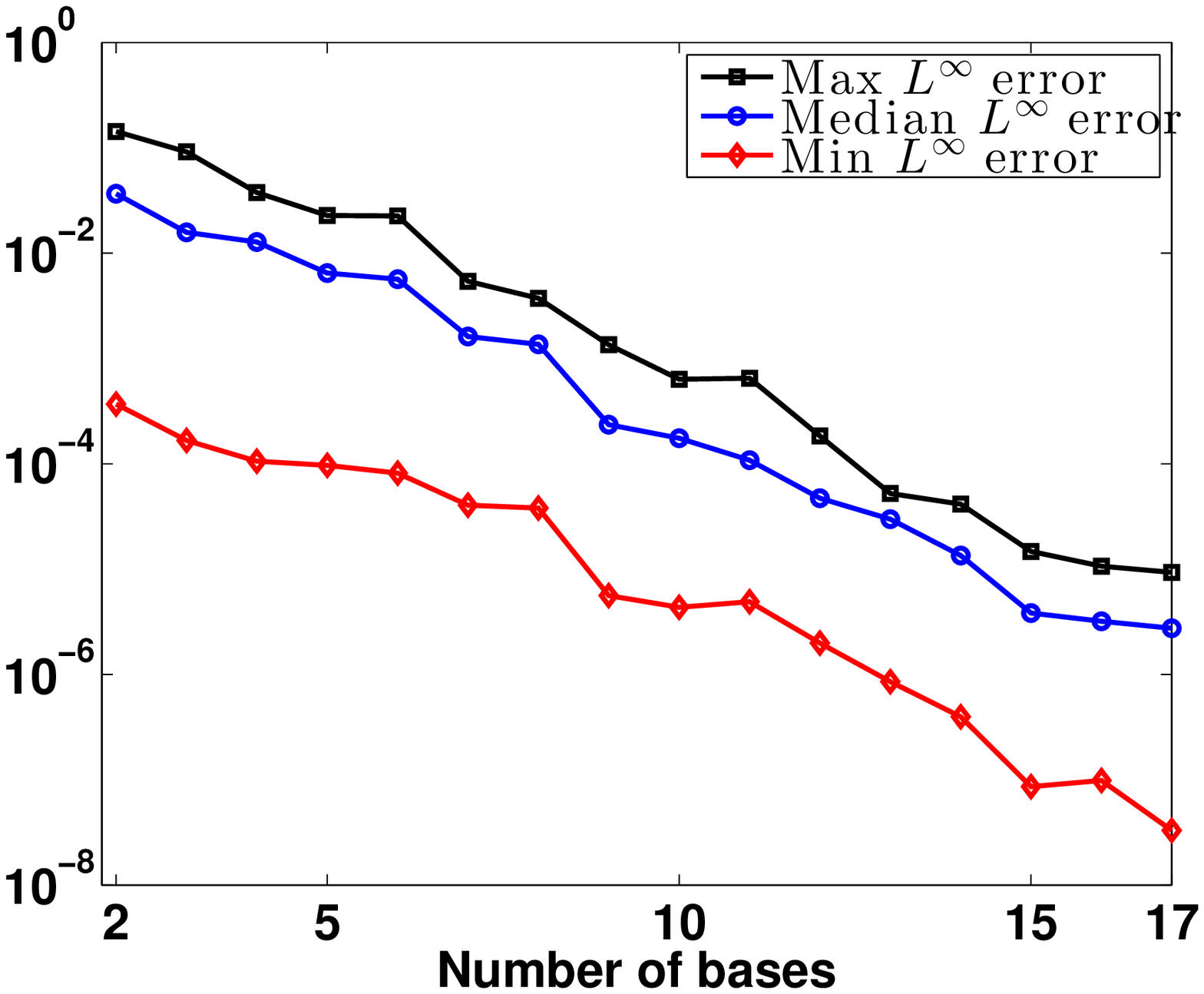}
        \includegraphics[width=0.45\textwidth]{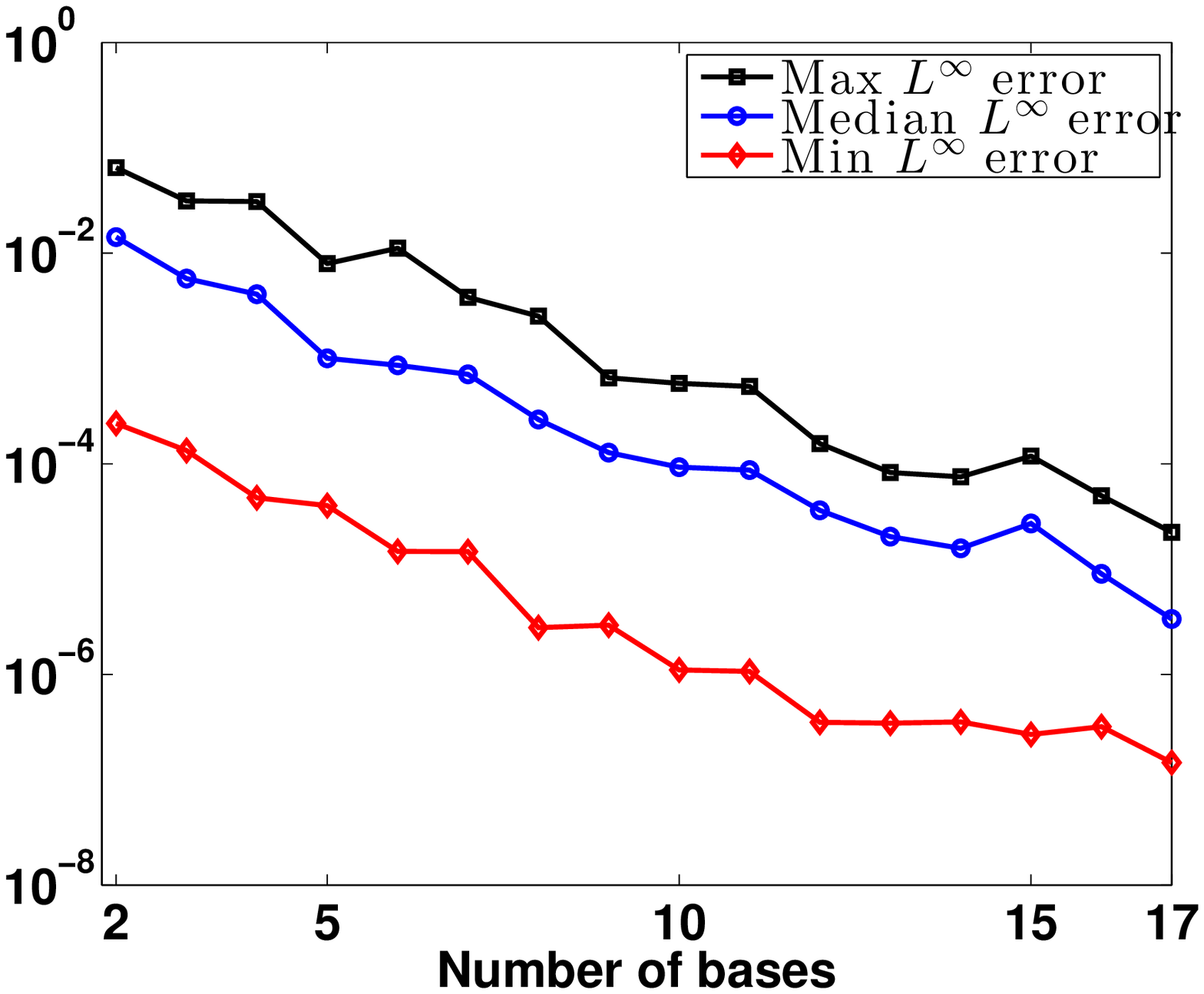}
  \end{center}
    \caption{History of convergence for the error estimate (top), the $L^2$ error (middle) and
    the $L^\infty$ error (bottom)of the RBM solutions for the anisotropic wavespeed simulation. On the left are for
    the least squares reduced collocation method, and the empirical reduced collocation results are on the right.}
    \label{fig:histconv}
\end{figure}

We first perform the offline pre-computation of the reduced basis  and collocation points.
The $17$ parameter values are chosen from $\Xi$   by Algorithms \ref{alg:LSgreedy} and \ref{alg:RCgreedy} are shown in Figure
\ref{fig:pickedmu},  with larger marker indicating the earlier that parameter picked. The reduced set of collocation
points $C_R^N$ for ERCM is shown in Figure \ref{fig:diffecmAniso} (top left). $C_R^i$ contains the $i$ points in the
computational domain $\Omega = [-1,1]^2$ corresponding to the $i$ largest markers.

Next, we solve for the reduced basis solution for a randomly selected set of $2,097$ parameter values in $\calD$
and compute the maximum, median, and minimum  errors for each selected value between the reduced solution and the
truth approximations. These, together with the maximum of the error estimate are plotted in Figure
\ref{fig:histconv}. We clearly see exponential convergence in all cases by both methods. We compare Figure
\ref{fig:truthaccuracy} and Figure \ref{fig:histconv} to draw the following remarkable conclusion: In the worst
case scenario, using the empirical reduced collocation method on a $16 = 4 \times 4$ grid can produce
solution having comparable accuracy of the truth approximation on a grid $50 \times 50$. We also see that, on average, the two proposed algorithms have similar accuracy. But, over a wide range
of parameter values, the least squares approach seems to be more stable (the errors have smaller variance).
Moreover, we show in Figure \ref{fig:diffecmAniso} how the choice of the reduced set of collocation points
$C_R^N$ affects the accuracy of the reduced collocation method: our proposed method generates the reduced
grid on the top left. The best case scenario for a randomly selected set of $2,097$ parameter values are
shown on the top right. On the other hand, if we naively use a coarse Chebyshev grid as the $C_R^N$ (shown
bottom left), the best case convergence plot is on the bottom right: the approximation is very bad with the
system becoming numerically singular for $N > 10$.

%\clearpage

\subsection{Results of the reduced solver: Diffusion}

We set $\calD = [-0.99,0.99]^2$, apply the empirical and least squares reduced collocation methods to the
diffusion problem and present the results in this section.

\begin{figure}[ht]
  \begin{center}
        \includegraphics[width=0.49\textwidth]{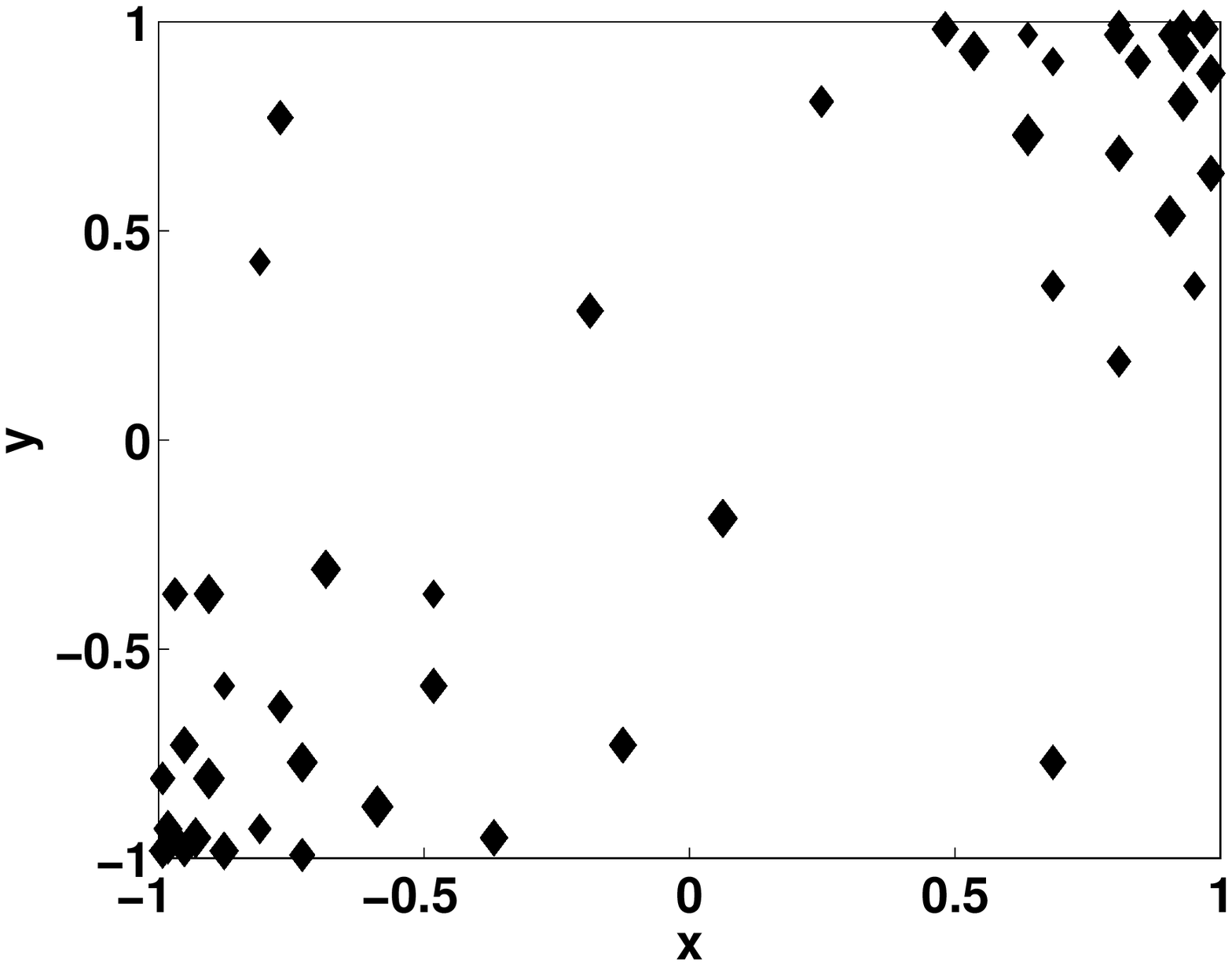}
  \end{center}
    \caption{The {\em reduced collocation points} selected by the ECM for diffusion problem.}
    \label{fig:diffecmEllip}
\end{figure}

\begin{figure}[ht]
  \begin{center}
        \includegraphics[width=0.49\textwidth]{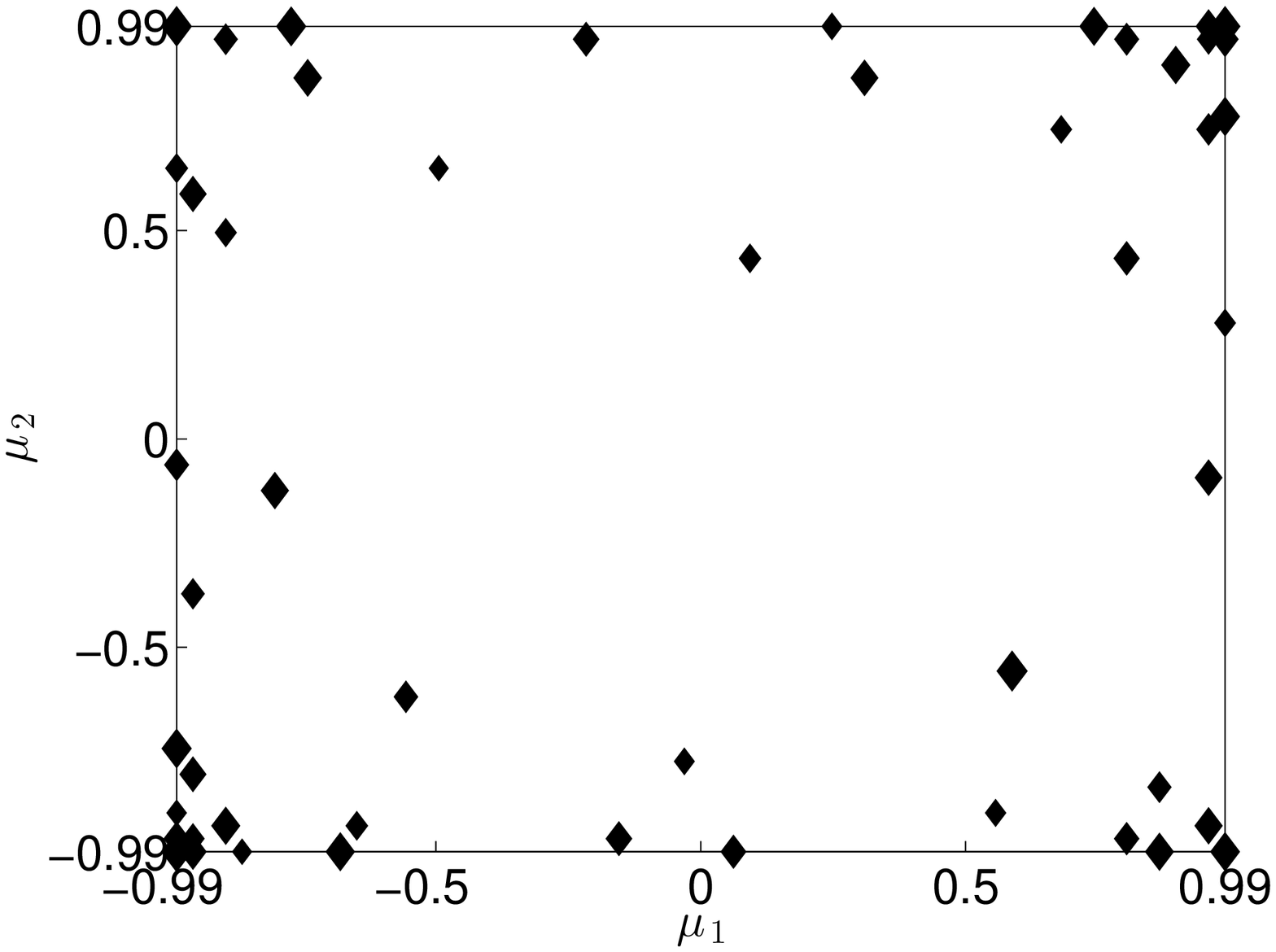}
        \includegraphics[width=0.49\textwidth]{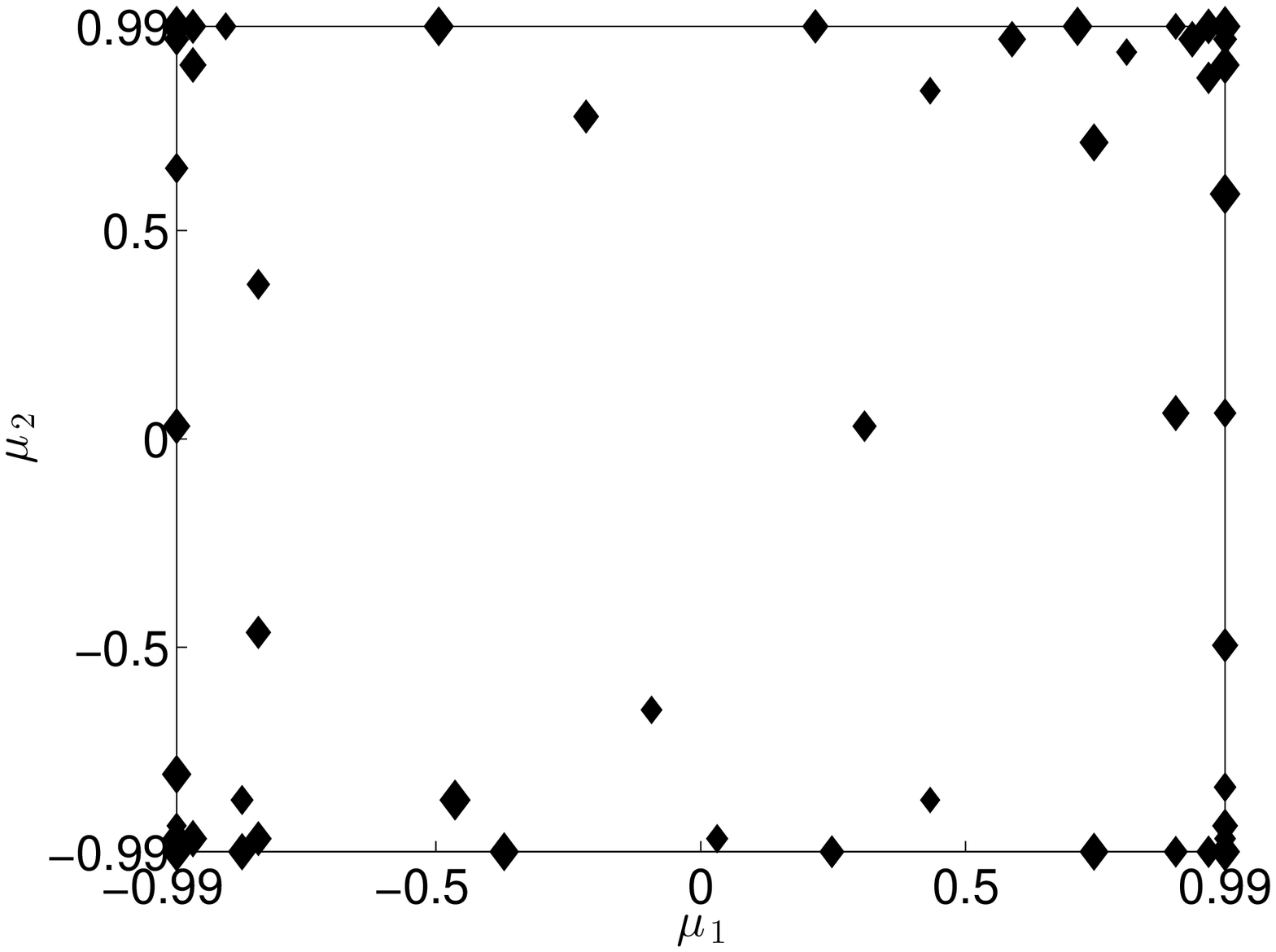}
  \end{center}
    \caption{The {\em parameter values} picked by the greedy algorithm for the diffusion equation. Left: Least squares. Right: reduced collocation.}
    \label{fig:diffpickedmu}
\end{figure}

\begin{figure}[ht]
  \begin{center}
        \includegraphics[width=0.45\textwidth]{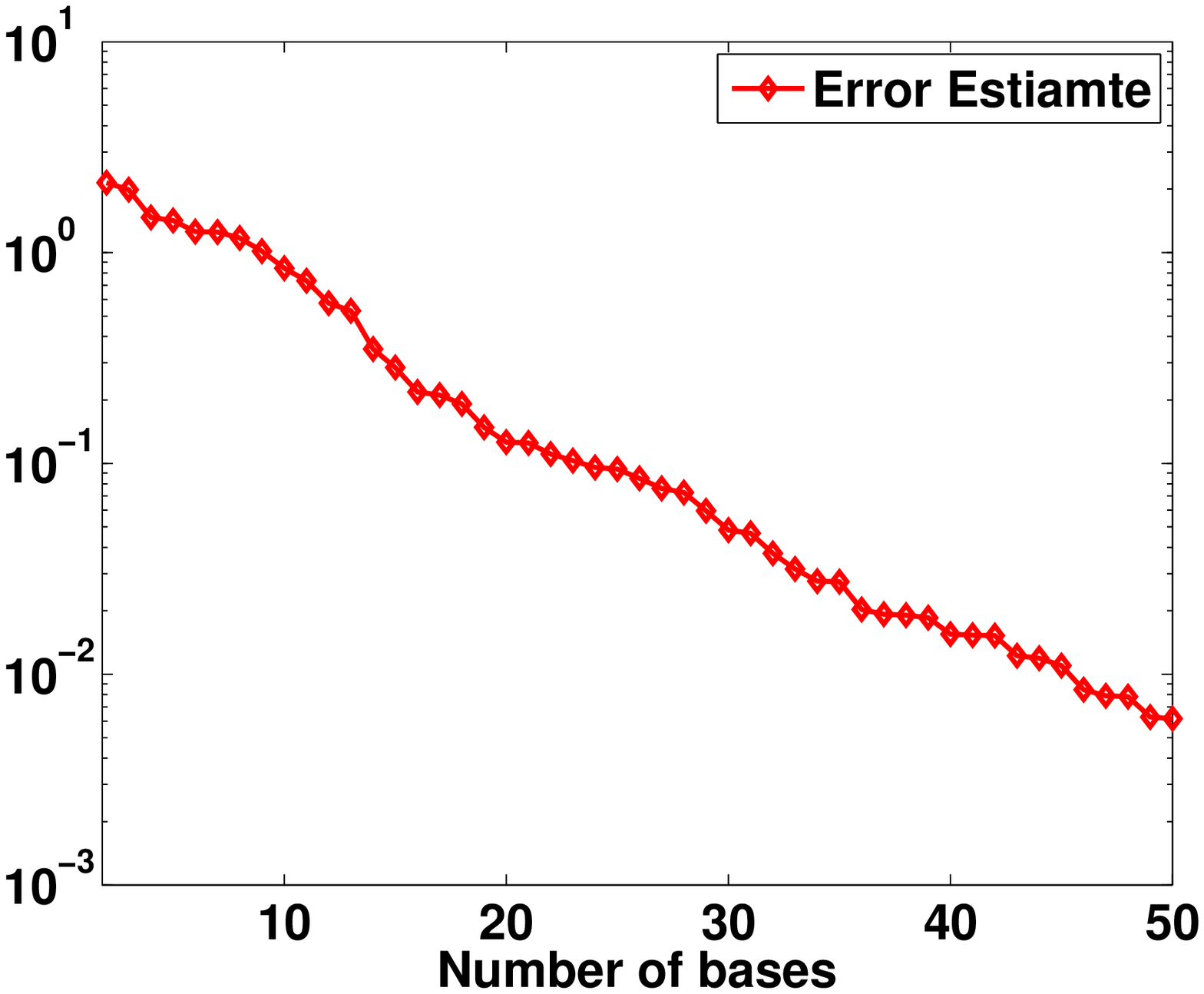}
        \includegraphics[width=0.45\textwidth]{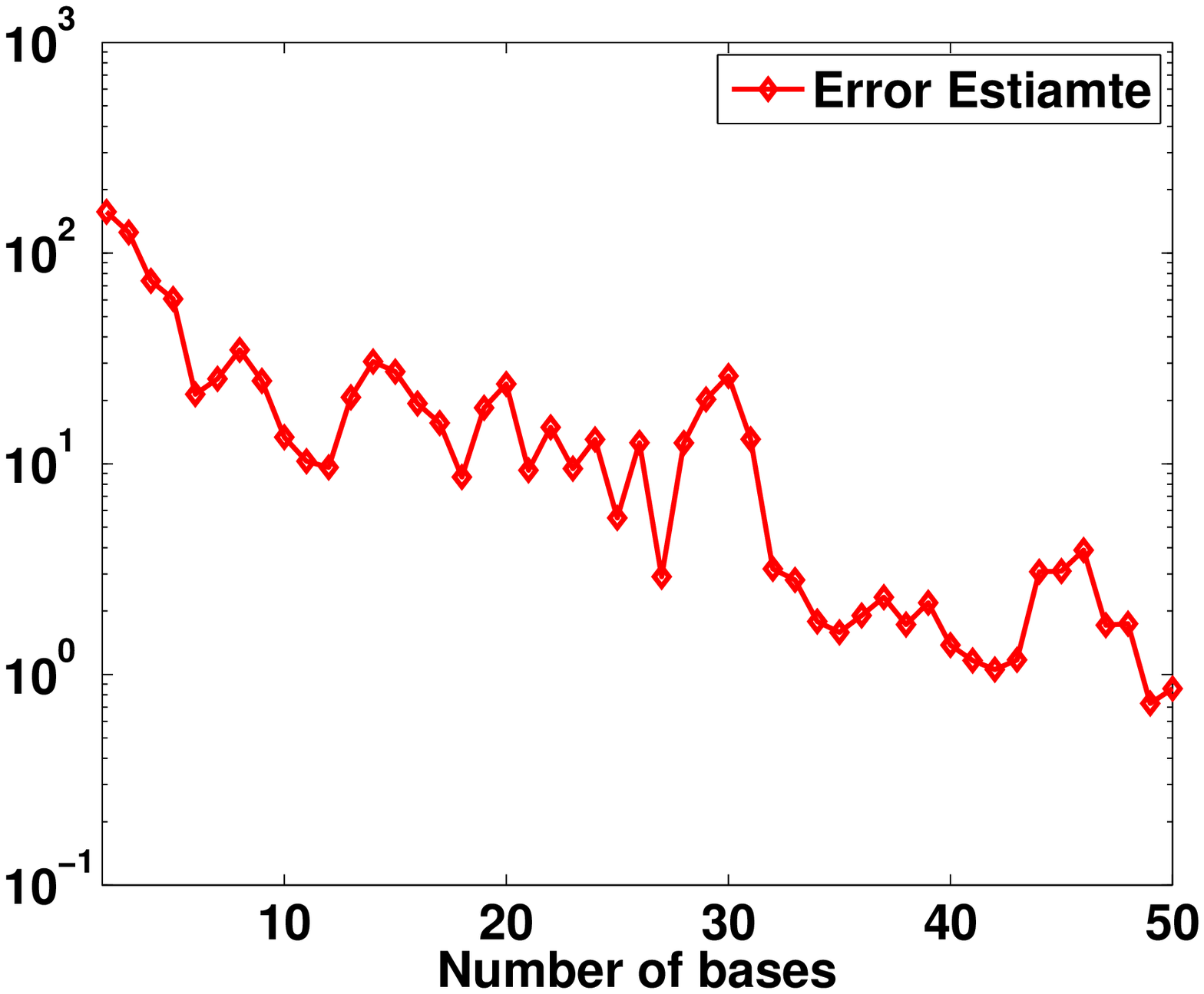}\\
        \includegraphics[width=0.45\textwidth]{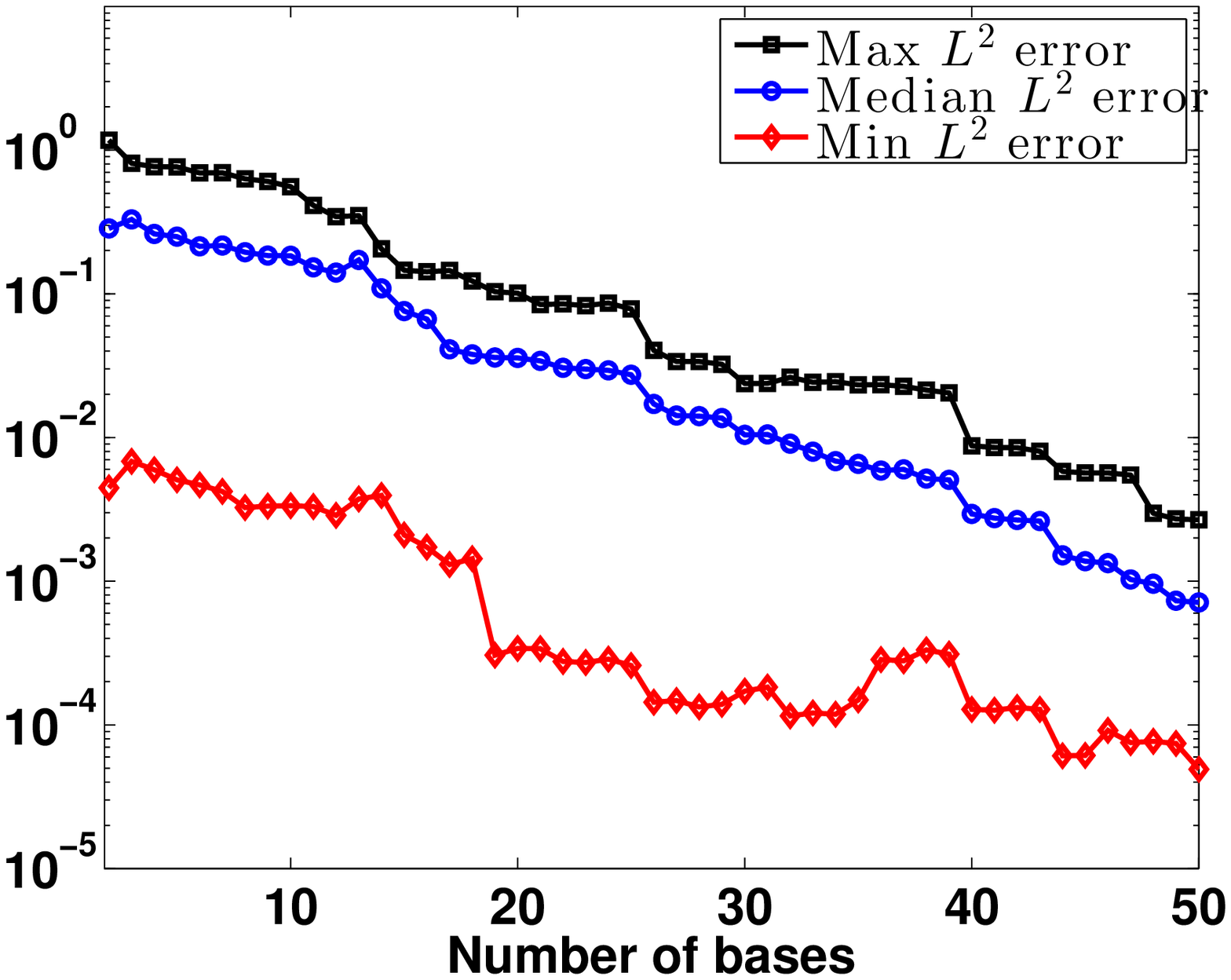}
        \includegraphics[width=0.45\textwidth]{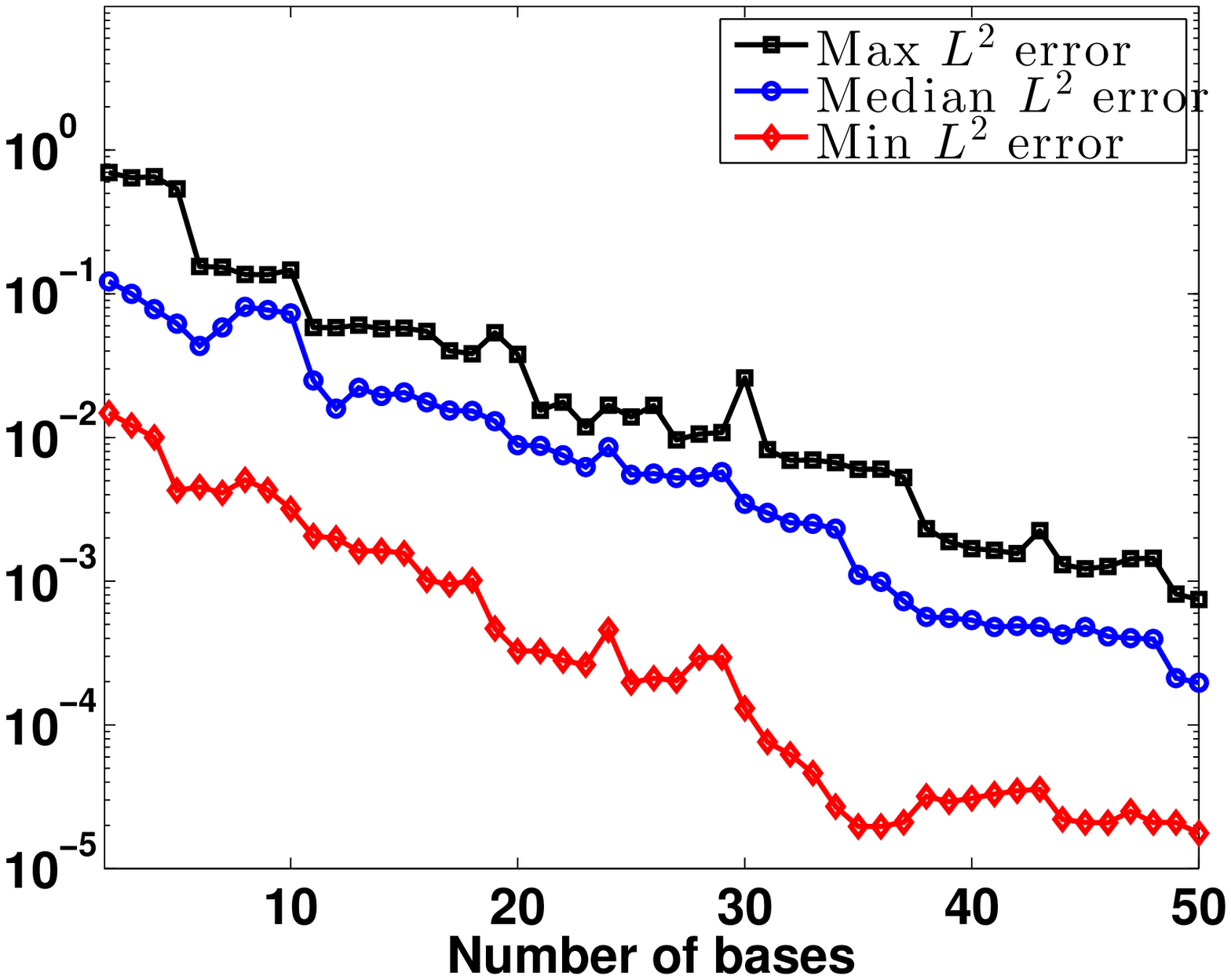}\\
        \includegraphics[width=0.45\textwidth]{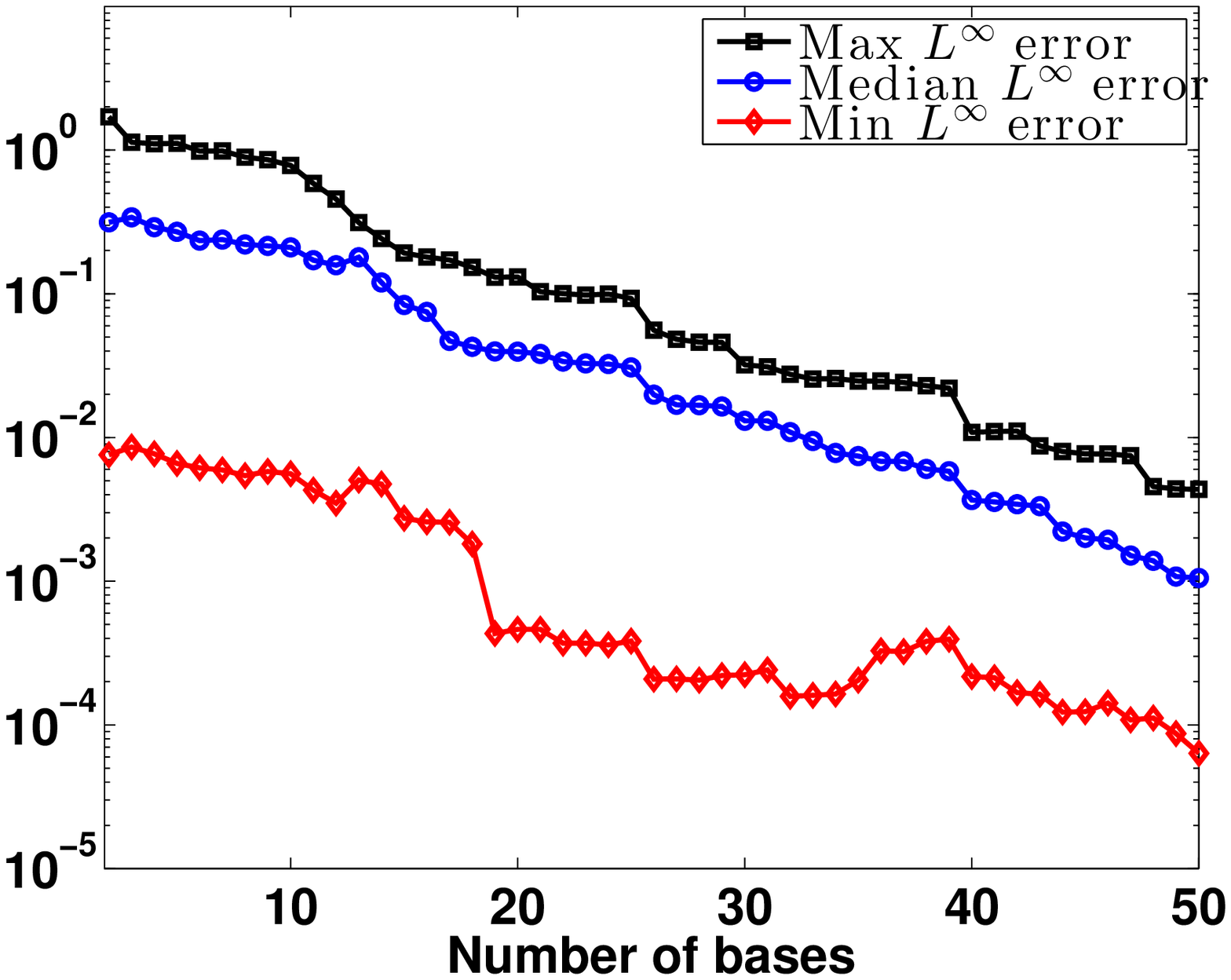}
        \includegraphics[width=0.45\textwidth]{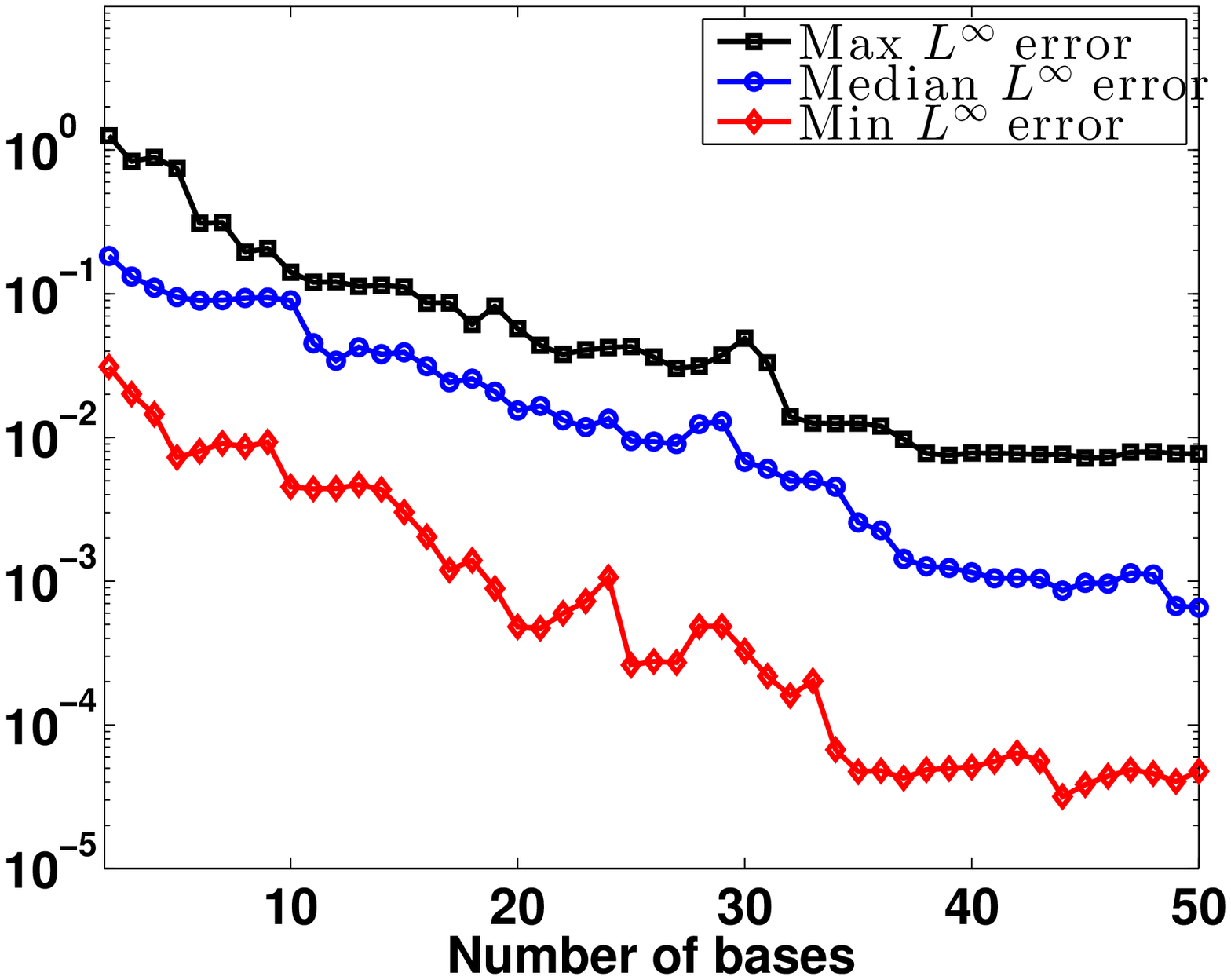}
  \end{center}
    \caption{History of convergence for the error estimate (top), the $\ell^2$ error (middle) and
    the $L^\infty$ error (bottom)of the RBM solutions for the diffusion problem. On the left are for
    the least squares approach, and the reduced collocation results are on the right.}
    \label{fig:diffhistconv}
\end{figure}

We pick $50$ parameter values in $\calD$ according to the greedy algorithm. The result is in Figure
\ref{fig:diffpickedmu} with larger marker indicating the earlier it is picked. Correspondingly, the $50$
points in $\Omega$ determined by the ERCM for empirical reduced collocation are shown in Figure
\ref{fig:diffecmEllip}.

Next, we solve for the RB solutions for randomly selected set of $1,057$ parameter values and compute the
maximum of the errors for all selected between the reduced solution and the truth approximations. These,
together with the maximum of the error estimate are plotted in Figure \ref{fig:diffhistconv}. We clearly see
exponential convergence in all cases for both methods.

\subsection{Computation time of the reduced solver}

In this section, we present statistics of the computation time for the reduced collocation methods. We
present in Table \ref{tab:comptime} the offline and online computational time. We normalize the time with
respect to that for solving truth approximation once. We see that the algorithm achieves savings of three
orders of magnitude. From these examples, it seems that the empirical collocation approach is a little
more efficient than the least-squares approach.

\begin{table}
\begin{center}
\begin{tabular}{|l|c|c|c|}
\hline Method & Offline time & Online time for $u_\mu^{(N)}$ & Time for $u_\mu^\calN$ \\
\hline Anisotropic LSRCM & $1.029 \, \times \, 10^3$ &  $2.18 \,  \times 10^{-4}$ & $1$\\
Anisotropic ERCM & $0.996 \,  \times 10^3$ & $2.62 \,  \times10^{-4}$ & $1$\\
Diffusion LSRCM & $8.265 \,  \times 10^3$ & $9.71 \,  \times10^{-4}$ & $1$\\
Diffusion ERCM & $7.730 \,  \times 10^3$ & $10.69 \,  \times10^{-4}$ & $1$\\
\hline
\end{tabular}
\end{center}
\caption{Computation times of the methods for the two test problems.} \label{tab:comptime}
\end{table}

\section{Concluding Remarks}
\label{sec:conclude}

In this paper, we propose {\em the first} reduced basis method for the collocation framework. Two rather
different approaches have been proposed and tested. They are both Galerkin-free but produce the same fast
exponential convergence and speedup as for the traditional Galerkin approach. In future work, we will examine
the accuracy and efficiency of our proposed methods for non-affine and nonlinear problems. We also plan to
study and tailor successive constraint method, currently used for computation of the lower bound for the
eigenvalues in the Galerkin setting \cite{HuynhSCM,CHMR-Cras,HKCHP}, for the collocation setting. \rev{It is also very interesting to apply the methods to more general collocation methods and to perform a detailed numerical comparison between the Galerkin approach and the collocation approaches introduced in this paper.}

\section*{Acknowledgements} The authors would like to thank
Professor Maday, Yvon from Paris VI University for helpful discussions that led to a deeper understanding of
the strength of our proposed approach. They also wish to thank the anonymous referees for constructive criticism that led
to an improved presentation of the material in this paper.

\end{document}